\documentclass[reqno]{article}
\usepackage{amsmath}
\usepackage{amsthm}
\usepackage{amssymb}
\usepackage{amscd}
\usepackage[all]{xy}
\usepackage{graphicx}

\newtheorem{thm}{Theorem}[section]
\newtheorem{prop}[thm]{Proposition}
\newtheorem{cor}[thm]{Corollary}
\newtheorem{lem}[thm]{Lemma}
\newtheorem{conj}[thm]{Conjecture}
\newtheorem*{clm}{Claim}

\theoremstyle{definition}
 \newtheorem{df}[thm]{Definition}
 
 \newtheorem{rem}[thm]{Remark}

 \newtheorem{fact}[thm]{Fact}
 \newtheorem{expl}[thm]{Example}

\newcommand{\tht}[2]{\vartheta
    \genfrac[]{0pt}{}{0}{#1} (k^{-1}\Omega, #2)}
\newcommand{\dgh}{d_{\mathrm{GH}}}
\newcommand{\xt}{{}^t}
\newcommand{\lt}{{\rule[-2.3ex]{0pt}{5.6ex}}^t \!\!\!}
\newcommand{\mt}{{\rule[-1.6ex]{0pt}{4.25ex}}^t \!\!\!}
\newcommand{\vin}{\cup \kern-0.8ex \hbox{\rule[0.03ex]{0.08ex}{1.35ex}}}

\newcommand{\dr}{\partial}
\newcommand{\db}{\bar{\partial}}
\newcommand{\dd}{\dr \db}
\newcommand{\diff}{\mathrm{Diff}(M)}
\newcommand{\symp}{\mathrm{Symp}(M,\omega)}

\newcommand{\lie}{\mathrm{Lie}\,}
\newcommand{\ep}{\varepsilon}
\newcommand{\cp}[1]{\mathbb{CP}^{#1}}
\newcommand{\fs}{\omega_{\mathrm{FS}}}
\newcommand{\im}{\mathrm{Im}\,}
\newcommand{\re}{\mathrm{Re}\,}

\begin{document}

\title{Lagrangian fibrations and theta functions}
\author{Yuichi Nohara\thanks{%
      e-mail : m98014i@math.nagoya-u.ac.jp}}
\date{Graduate School of Mathematics,\\ Nagoya University}

\maketitle

\begin{abstract}
In this thesis we study asymptotic behavior of projective embeddings
of abelian varieties and their amoebas.
The projective embeddings are given by theta functions.
It is known that a Lagrangian fibration of the abelian variety 
determines a basis of theta functions.
After reviewing the relation from the viewpoint of geometric quantization
and mirror symmetry,
we prove that the Lagrangian fibration of the abelian variety can be approximated
by the moment maps of natural torus action in a suitable sense.
In the final section we discuss a part of the result in more general setting.
\end{abstract}

\tableofcontents

\newpage
\section{Introduction}

The purpose of this thesis is to discuss some relations between
theta functions and Lagrangian fibrations of abelian varieties.

Lagrangian fibration is a map $(M,\omega) \to B$ from a symplectic 
manifold such that each (smooth) fiber is a Lagrangian submanifold.
Lagrangian fibrations appear in several areas of mathematics and
mathematical physics such as completely integrable systems 
(classical mechanics), geometric quantization and mirror symmetry.
In particular, the picture of mirror symmetry via 
{\it special Lagrangian fibrations}
proposed by Strominger-Yau-Zaslow \cite{SYZ} has attracted much attention
and a lot of works relating to this program
have been carried on
(e.g. \cite{G1,G2}, \cite{Hi}, \cite{R1,R2,R3}).

However, little is known about (special) Lagrangian submanifolds and
fibrations at the moment.
In general, it is more difficult to deal with submanifolds compared
to functions or connections.
In the case of (special) Lagrangian fibrations,
the main difficulties seems to be the following two points.
The first one is the fact that 
general Lagrangian fibrations have singular fibers.
Singularities appearing in Lagrangian fibrations are not 
completely understood yet.
The other one concerns the case of special Lagrangian fibrations.
To define the special Lagrangian submanifolds, 
we need a Ricci-flat metric.
Ricci-flat metrics are difficult objects to analyze.

The cases of Abelian varieties and toric varieties are 
typical examples of Lagrangian fibrations which are well understood.
Mirror symmetry via special Lagrangian fibrations
works well for abelian varieties and theta functions 
play an important role in the theory(\cite{PZ}, \cite{F}).

Theta functions, or more generally holomorphic sections of 
ample line bundles $L$ on a smooth projective variety $X$,
are the other main characters in this paper.
Of course, the space $H^0(X,L^k)$ of holomorphic sections
does not have natural basis in general.
However, some projective varieties have natural basis.
For example, in the case of toric varieties, monomials are natural
basis of holomorphic sections.
Also in the case of abelian varieties, there are natural basis of
theta functions(one can find such basis in the Mumford's book \cite{M1}).
We can think that
such basis are determined by Lagrangian fibrations
by using the notion of geometric quantization
(the case of abelian varieties is discussed in \cite{FMN},
\cite{Ty2}, \cite{W} and other cases are, for example, in 
\cite{A}, \cite{BPU}, \cite{GS}, \cite{H}, \cite{JW}).
The case of abelian varieties can be also interpreted 
as a special case of mirror symmetry.

Note that, for large $k$, each basis of $H^0(X,L^k)$ defines an embedding
into the projective space $\cp{N_k} = \mathbb{P} H^0(X,L^k)^*$:
\[
  \iota_k : X \hookrightarrow \cp{N_k},
  \quad z \mapsto (s_0(z): \dots : s_{N_k}(z)) 
\]
We see the relation between Lagrangian fibration of $X$ and
projective embeddings defined by the basis corresponding to the 
Lagrangian fibration.
Namely, we compare the Lagrangian fibration with
the standard one of projective spaces (i.e. the moment maps of natural
torus actions)
\[
  \mu_k : \cp{N_k} \longrightarrow \mathbb{R}^{N_k} 
  = \mathrm{Lie}\,T^{N_k}.
\]
We denote the restriction of the moment map to $X$ by
\[
  \pi_k = \mu_k \circ \iota_k : X \longrightarrow B_k \, ,
\]
where $B_k$ is the image of $X$ under the moment map.

The simplest example is the case of toric varieties.
If we take monomials as a basis of $H^0(X,L^k)$,
the corresponding embeddings are torus equivariant.
This means that $\pi_k : X \to B_k$ coincides with
the moment map of $X$.

The main theorem of this paper deal with the case of abelian varieties.
In this case, the situation is not so trivial.
The rough statement is the following.

\begin{thm} 
  Let $X$ be an abelian variety. 
  We consider projective embeddings oof $X$
  by using theta functions determined by a Lagrangian fibration
  of $X$.
  Then the sequence $\{ \pi_k : X \to B_k \}$ converges to the original
  Lagrangian fibration as $k \to \infty$
  in the ``Gromov-Hausdorff topology''.
\end{thm}
The precise statement is given in Section 6.

This theorem is closely related to the theory of approximation of
K\"ahler metrics by Fubini-Study metrics.
Tian \cite{T} and Zelditch \cite{Zel} proved that
every K\"ahler metric in a fixed K\"ahler class $c_1(L)$ can be 
approximated by restrictions of (normalized) Fubini-Study metrics
\[
  \omega_k = \frac 1k \iota_k^* \fs
\]
on $\cp{N_k} = \mathbb{P} H^0(X,L^k)^*$ for large $k$
by choosing appropriate embeddings.

There is another theory on this subject and it is related to the notion
of stability in the sense of geometric invariant theory.
This is a possible approach to prove ``Hitchin-Kobayashi correspondence for
manifolds'' which claims the equivalence of stability of a projective
variety and the existence of special K\"ahler metrics such as 
K\"ahler metric of constant scalar curvature.
In this approach, we relate the choice of projective embeddings with
the stability condition and study the asymptotic behavior of 
the restrictions of Fubini-Study metrics.
In particular, this theory could be effective for the the study of
Ricci-flat metrics.

The projective embeddings defined by theta functions 
give examples of both of the above theories.
The main theorem can be regarded as a prototype of approximations of
K\"ahler metrics and Lagrangian fibrations at the same time.
Our method might be useful for the study of special Lagrangian fibrations.

This paper is organized as follows.
In Section 2, we recall basic facts on symplectic geometry used
in this paper.
In Section 3, we review the theory of geometric quantization,
especially, the part related to Lagrangian fibration in detail.
Section 4 is devoted to a brief summary of mirror symmetry for 
abelian varieties via special Lagrangian fibrations.
We recall the theorem of Tian and Zelditch stated above in Section 5.
In this section, we also recall the relation between projective 
embeddings, stability of polarized manifolds and canonical metrics.
The precise statement and proof of the main theorem is given in 
Section 6.
Finally we discuss a generalization of the above result in Section 7.

\paragraph{Acknowledgment}
The author would like to thank Ryoichi Kobayashi for helpful discussions.
A large part of this thesis is based on the talk at
``NIT Homotopy Theory Meeting''.
The author also thanks the organizer Norihiko Minami.


\newpage
\section{Basic facts on symplectic geometry}

In this section, we summarize basic facts on symplectic geometry
which we need later.


\subsection{Symplectic manifolds and automorphisms}

\begin{df}
  Let $M$ be a smooth manifold of dimension $2n$.
  A {\it symplectic structure} $\omega$ on $M$ is a
  non-degenerate closed 2-form on $M$.
  A {\it symplectic manifold} is a pair $(M,\omega)$ of
  smooth manifold and symplectic structure on it.
\end{df}

We denote the group of diffeomorphisms of $M$ by $\diff$.
Its Lie algebra $\lie \diff$ is identified formally with the Lie algebra
$\mathcal{X} (M)$ of vector fields on $M$.
By using the symplectic structure $\omega$, we can identify
$\mathcal{X} (M)$ with the space of 1-forms on $M$:
\[
  \mathcal{X} (M) \longrightarrow \Omega^1(M) \,, \quad
  \xi \longmapsto i_{\xi} \omega \, .
\]
Let $\symp \subset \diff$ be the subgroup of diffeomorphisms
which preserve the symplectic structure $\omega$.
Diffeomorphisms in $\symp$ are called {\it symplectomorphisms}.
\begin{prop}
  Under the above identification, the Lie algebra of $\symp$
  corresponds to the space $Z^1(M)$ of closed 1-forms.
\end{prop}

This proposition follows from $d \omega = 0$ and the following formula
\[
  L_{\xi} \omega = i_{\xi} d \omega + d (i_{\xi} \omega)\, ,
\]
where $L_{\xi}$ is the Lie derivative.

\begin{df}
  For a smooth function $f \in C^{\infty}(M)$,
  its {\it Hamilton vector field} $\xi_f$ is defined by
  \[
    i_{\xi_f} \omega = df \, .
  \]
  A {\it Hamilton diffeomorphism} is a diffeomorphism generated by
  Hamilton vector fields.
  We denote the group of Hamilton diffeomorphisms by
  $\mathrm{Ham}\,(M,\omega)$.
\end{df}
From the definition,
the Lie algebra of $\mathrm{Ham}\,(M,\omega)$ is identified with
the space $B^1(M)$ of exact 1-forms.
In particular, every Hamilton diffeomorphism preserves $\omega$.

\begin{prop}
  Assume that $\omega$ represents an integral cohomology class:
  \[
     [\omega] \in H^2(M, \mathbb Z)\, .
  \]
  Then there exists a Hermitian line bundle $L \to M$
  with a unitary connection $\nabla$ such that
  $c_1(L,\nabla) = \omega$.
\end{prop}
Such $(L, \nabla)$ is called a {\it prequantum bundle} of $(M,\omega)$.
This is a fundamental ingredient in geometric quantization.

\begin{proof}
  Let $\mathcal{A}$ (resp. $\mathcal{A}^*$) be the sheaf of smooth
  functions (resp. non-zero smooth functions) and consider the
  following exact sequence:
  \[
   \begin{CD}
    0 @>>> \mathbb{Z} @>>> \mathcal{A}
      @>{\exp(2\pi \sqrt{-1} \cdot)}>>
      \mathcal{A}^* @>>> 0\, .
   \end{CD}
  \]
  Recall that $H^1(M, \mathcal{A}^*)$ parameterizes the isomorphism
  classes of smooth complex line bundles on $M$.
  Since $\mathcal{A}$ is a fine sheaf, the long exact sequence becomes
  \[
    \dots \longrightarrow H^1(M, \mathcal{A})
    \longrightarrow H^1(M, \mathcal{A}^*)
    \overset{c_1}{\longrightarrow}
    H^2(M, \mathbb{Z}) \longrightarrow 0\, .
  \]
  This implies that there exists a line bundle $L \to M$ such that
  $c_1(L) = [\omega]$.

  Now we take an arbitrary Hermitian metric and a unitary connection 
  $\nabla'$ on $L$.
  Note that its curvature $\omega'$ represents $c_1(L)$.
  There exists a smooth 1-form $\alpha$ such that
  \[
    \omega = \omega' + d \alpha\, .
  \]
  Then $\nabla = \nabla' - 2 \pi \sqrt{-1} \alpha$ 
  gives a unitary connection
  such that its curvature coincides $\omega$.
\end{proof}

We consider the automorphism group $\mathcal{G}$ of
a prequantum bundle $(L,\nabla) \to (M,\omega)$.
$\mathcal{G}$ consists of isomorphisms
\begin{equation*}
 \begin{matrix}
  L & \stackrel{\hat{F}}{\longrightarrow} & L \\
  \downarrow & & \downarrow \\
  M & \stackrel{F}{\longrightarrow} & M
 \end{matrix}
\end{equation*}
such that $\hat{F}$ preserves the Hermitian metric and the connection
$\nabla$ of $L$ (hence $F$ preserves $\omega$).

\begin{prop} \label{lie}
  The Lie algebra of $\mathcal{G}$ is formally identified with the space of
  smooth functions $C^{\infty}(M)$ on $M$, here the Lie algebra structure
  of $C^{\infty}(M)$ is given by the Poisson bracket.
\end{prop}

We can see this as follows.
Consider the natural projection
\[
  \mathcal{G} \longrightarrow \mathrm{Symp}\,(M,\omega)\, ,
  \quad (\hat{F},F) \longmapsto F
\]
and denote its image by $\mathcal{G}_0$.
The kernel consists of automorphisms which preserve the base space.
It is easy to see that this is isomorphic to $S^1$.
Then we have the following exact sequence:
\[
  \begin{CD}
    1 @>>> S^1 @>>> \mathcal{G} @>>> \mathcal{G}_0 @>>> 1 \, .
  \end{CD}
\]
Proposition \ref{lie} is a consequence of the fact that 
$\mathcal{G}_0 \cong \mathrm{Ham}\,(M,\omega)$.
In fact, the corresponding exact sequence of Lie algebras
is given by
\[
  \begin{CD}
    0 @>>> \mathbb{R} @>>> \mathrm{Lie}\,\mathcal{G} @>>>
    d \, C^{\infty}(M) @>>> 0 \, .
  \end{CD}
\]

Let $\Gamma (L) = \Gamma (M, L)$ be the space of smooth sections
of $L$.
Then,
for each $f \in C^{\infty}(M) \cong \mathrm{Lie} \mathcal{G}$, its action 
$ \hat{f} : \Gamma (L) \to \Gamma (L)$ is given by
\begin{equation}
  \hat{f} (s) = \nabla_{\xi_f} s - \sqrt{-1} fs \, .
  \label{f-hat}
\end{equation}


\subsection{Action-angle coordinate}

In this and the next subsection, 
we recall basic properties of Lagrangian fibrations.
\begin{df}
 A submanifold $S \subset M$ is called a
 {\it Lagrangian submanifold} if it satisfies $\omega |_S =0$
 and $\dim S = \frac12 \dim M$.
\end{df}

We call a fibration $\pi : M \to B$ a {\it Lagrangian fibration}
if general fibers are Lagrangian.

\begin{rem}
It is natural to allow Lagrangian fibrations to have degenerate
fibers.
However we do not care about singular fibers in this section.
\end{rem}

\begin{expl}
  Let $B$ be an $n$-dimensional manifold and $T^*B$  its
  cotangent bundle.
  Then the natural projection $T^*B \to B$ is a Lagrangian fibration
  with respect to the standard symplectic structure $\omega_0$
  on $T^*B$.
\end{expl}

Now we recall the notion of moment maps.
\begin{df}
Let $G$ be a Lie group acting on a symplectic manifold
$(M,\omega)$ as symplectomophisms.
A {\it moment map} 
\[
  \mu : M \longrightarrow \mathfrak{g}^* = (\mathrm{Lie}\, G)^*
\]
is a $G$-equivariant map satisfying
\[
  d \langle \mu , \xi \rangle = i_{\xi}\omega
\]
for $\xi \in \mathfrak{g}$,
where $\langle \, , \, \rangle$ is the natural pairing and
we identified $\xi \in \mathfrak{g}$ with the vector field
on $M$ defined by its action.
\end{df}

Remark that moment maps 
\[
  \mu : M \longrightarrow \mathfrak{t}^* \cong \mathbb{R}^n
\]
of $T^n$-action ($n = \frac 12 \dim M$) are Lagrangian fibrations
(with degenerate fibers in general).

\begin{prop} \label{action-angle}
  Let $\pi : (M,\omega) \to B$ be a Lagrangian fibration.
  \begin{enumerate}
    \item Every smooth fiber has a natural affine structure.
          In particular, compact fiber is a torus $T^n$.
    \item Assume that $\pi^{-1}(b_0)$ is smooth and compact.
          Then there exist a neighborhood $U \subset B$ of $b_0$
          and local coordinates $(x^1, \dots , x^n)$ of $U$
          and $(y^1, \dots , y^n)$ of $T^n$ such that
      \begin{enumerate}
        \item $(x^1, \dots , x^n,y^1, \dots , y^n)$ is a coordinate
              of $\pi^{-1}(U) \cong U \times T^n$ with
              $\omega = \sum dx^i \wedge dy^i$, and
        \item there exists a $T^n$-action on $\pi^{-1}(U)$
              such that
              \[
                (x^1, \dots , x^n) : \pi^{-1}(U) \to \mathbb{R}^n
              \]
              gives its moment map.
      \end{enumerate}
  \end{enumerate}
\end{prop}

Such coordinates $(x^1, \dots , x^n,y^1, \dots , y^n)$ are called
{\it action-angle coordinates}.

\paragraph{Outline of the proof.}
1. Let $U$ be a neighborhood of $b_0 \in B$.
We may assume that $U \subset \mathbb{R}^n$
is a small ball.
Denote the coordinate of $\mathbb{R}^n$ by
$(\tilde{x}^1, \dots , \tilde{x}^n)$ and set
$f_i = \tilde{x}^i \circ \pi : \pi^{-1}(U) \to \mathbb{R}$.

\begin{clm}
  $\{f_i, f_j\} = 0$ for each $i,j = 1 \dots , n$.
\end{clm}

\begin{proof}
Since $f_i$ is constant along fibers,
\[
  0 = \xi f_i = df_i (\xi) = \omega (\xi_{f_i}, \xi)
\]
for each $\xi \in T \bigl( \pi^{-1}(b) \bigr)$.
In particular, $\xi_{f_i} \in T \bigl( \pi^{-1}(b) \bigr)$ and
\[
  0 = \omega (\xi_{f_i}, \xi_{f_j}) = \{f_i, f_j\} \, .
\]
\end{proof}

Note that $\xi_{f_1}, \dots , \xi_{f_n}$ are linearly independent.
Furthermore, we have
\[
  [ \xi_{f_i}, \xi_{f_j} ] = \xi_{\{f_i, f_j\}} = 0
\]
from the above claim.
Then $\xi_{f_1}, \dots , \xi_{f_n}$ generates an
$\mathbb{R}^n$-action on $\pi^{-1}(U)$ and
each fiber $\pi^{-1}(b)$ is preserved by this action.
In particular, if fibers are compact, we have $\pi^{-1}(b) \cong T^n$.

2. Here we only recall the construction of action-angle coordinates 
$(x^1, \dots , x^n,y^1, \dots , y^n)$.
Since $U$ is contractible, the inclusion
$\iota : \pi^{-1}(b_0) \hookrightarrow \pi^{-1}(U)$
induces an isomorphism
\[
  \iota^* : H^2(\pi^{-1}(U), \mathbb{R})
  \overset{\sim}{\longrightarrow}
  H^2(\pi^{-1}(b_0), \mathbb{R}) \, .
\]
Because $\pi^{-1}(b_0)$ is Lagrangian, we have
\[
  \iota^* [\omega] =
  \bigl[\omega |_{\pi^{-1}(b_0)} \bigr] = 0 \, .
\]
This implies that there exists a $1$-form $\theta$ on
$\pi^{-1}(U)$ such that $d \theta = \omega$.
We choose generators $\gamma_i^b$ of $H_1(\pi^{-1}(b),\mathbb{Z})$,
$i = 1, \dots ,n$, so that they depend smoothly on $b \in U$
and set
\[
  x^i(b) = \frac 1{2 \pi} \int_{\gamma_i^b} \theta \, .
\]
Note that $x^i(b)$ depends only on the homology classes
$\bigl[\gamma_i^{b_0} \bigr] \in H_1(\pi^{-1}(b_0), \mathbb{Z})$.

Then we can show that  $(x^1, \dots , x^n)$ defines a coordinate
around $b_0$.
Moreover, 
  $\exp (\xi_{x^i} )(x) = x$ for $x \in \pi^{-1}(U)$.
In other words, $\exp (t \xi_{x^i} )$ has period $1$.

Fix a Lagrangian section $\lambda : U \to \pi^{-1}(U)$
and define a map
\[
  \psi : T^*U \to \pi^{-1}(U)
\]
by $\psi (df(b)) = \exp (\xi_f ) \cdot \lambda (b)$.
We define the coordinate of
$ \exp (\xi_{\sum y^i \xi_{x^i}}) \cdot \lambda$ by
$(y^1, \dots , y^n)$.
Note that the zero section of $T^*U$ corresponds to $\lambda$.
Then the coordinate $(y^1, \dots, y^n)$
has period $1$.
\qed

\begin{rem} \label{holonomy}
  If there exists a quantum line bundle $(L,\nabla)$,
  we can take $\theta$ to be a connection 1-form of $\nabla$.
  In this case, $\exp (2 \pi \sqrt{-1} x^i)$ is the holonomy of
  $(L,\nabla)$ on each fiber.
\end{rem}


\subsection{Global structure of Lagrangian torus fibrations}

Next we discuss global structures of Lagrangian fibrations.
Here we assume that the Lagrangian fibration $\pi : M \to B$
has no degenerate fiber.
In the proof of theorem \ref{action-angle}, we constructed a map
$\psi : T^*U \to \pi^{-1}(U)$.
Define $\Lambda_U \subset T^*U$ by
\[
  \Lambda_U = \{ \alpha \in T^*U \, | \, \psi (\alpha) = \lambda \,\}
\]
Then this is a $\mathbb{Z}^n$-bundle spanned by
$dx^1, \dots , dx^n$.
For each $b \in U$, $\Lambda_b$ can be identified with
$H_1(\pi^{-1}(b), \mathbb{Z}) \cong
H^{n-1}(\pi^{-1}(b), \mathbb{Z})^*$.
In other words, we have
\[
  \Lambda_U \cong \bigl( R^{n-1}\pi_* \mathbb{Z} \bigr)^*.
\]
This means that $\Lambda_U$ defines a global $\mathbb{Z}^n$-bundle
$\Lambda \subset T^*B$.
Since $\Lambda$ is spanned locally by $dx^1, \dots , dx^n$,
$\Lambda$ is a Lagrangian submanifold in $(T^*B, \omega_0)$.
Hence $\omega_0$ induces a symplectic structure on $T^*B/\Lambda$
and the natural projection $T^*B/\Lambda \to B$ is 
a Lagrangian torus fibration.
We remark that $T^*B/\Lambda$ has a Lagrangian section 
(i.e. the zero section).
This is called the {\it Jacobian fibration} of the Lagrangian fibration
$\pi : M \to B$ (\cite{G1}).

Now assume that the Lagrangian fibration
$\pi : (M,\omega) \to B$ has a Lagrangian section
$\lambda : B \to M$.
Then the map $\psi$ extends to a global map
\[
  0 \longrightarrow \Lambda \longrightarrow T^*B
    \overset{\psi}{\longrightarrow} M \longrightarrow 0 \, .
\]
Hence we have an isomorphism $T^*B/\Lambda \to M$ of Lagrangian fibrations.

Next we consider the case that $\pi : M \to B$ does not necessarily
have a global Lagrangian section.
Take a open cover $M = \bigcup_i U_i$ of $M$ such that a Lagrangian section
$\lambda_i : U_i \to \pi^{-1}(U_i)$ and
an action-angle coordinate are defined on each $U_i$.
Then we have
\[
  0 \longrightarrow \Lambda|_{U_i} \longrightarrow T^*U_i
    \longrightarrow \pi^{-1}(U_i) \longrightarrow 0
\]
on each $U_i$ (i.e. $M \to B$ is locally isomorphic to
$T^*B/\Lambda \to B$).
$M \to B$ can be reconstructed by gluing $T^*U_i/\Lambda|_{U_i} \to U_i$.
The gluing is determined by Lagrangian sections
$\mu_{ij} : U_i \cap U_j \to T^*B/\Lambda$  by
\[
  \lambda_i (b) = \mu_{ij} (b) \cdot \lambda_j (b)
\]
for each $b \in U_i \cap U_j$.
Let $\mathcal{L} (T^*B/\Lambda)$ be the sheaf of Lagrangian sections
of $T^*B/\Lambda$.
Then $\{ \mu_{ij} \}$ defines a cohomology class
$[ \mu ] \in H^1(B, \mathcal{L} (T^*B/\Lambda))$.

We summarize as follows:

\begin{thm}[Duistermaat \cite{Du}]
  \begin{enumerate}
    \item If a Lagrangian fibration $M \to B$ has a Lagrangian
          section, then
          $M \to B$ is isomorphic to $T^*B/\Lambda \to B$
          as Lagrangian fibrations.
    \item The isomorphism classes of Lagrangian fibrations
          $M \to B$ with Jacobian fibration
          $T^*B/\Lambda \to B$ are parameterized by
          $H^1(B, \mathcal{L} (T^*B/\Lambda))$.
  \end{enumerate}
\end{thm}


\newpage
\section{Geometric Quantization}

Geometric quantization is a method to
construct representations geometrically.
In the case of abelian varieties,
we obtain the representations of finite Heisenberg groups.
This view point gives an interpretation of the relation 
between Lagrangian fibrations of abelian varieties and theta functions.

\subsection{Brief review of quantum mechanics}

First we recall a formulation of
quantum mechanics on Euclidean spaces.

Let $M = T^* \mathbb{R}^n = \mathbb{R}^{2n}$ be a
symplectic vector space with the standard symplectic structure
$\omega = \frac 1{\hbar} \sum_{i=1}^n dq^i \wedge dp_i$,
where $(q^1, \dots , q^n)$ is a coordinate on $\mathbb{R}^n$,
$(p_1, \dots , p_n)$ the canonical coordinate on the cotangent spaces
and $\hbar$ is the Planck's constant.

The goal is to correspond a representation of a suitable Lie subalgebra
of $C^{\infty}(M)$.
For each function $f$, we denote the corresponding operator by
$\hat{f} : \mathcal{H} \to \mathcal{H}$, where
$\mathcal{H}$ is a Hilbert space.
From the requirement that $f \mapsto \hat{f}$ is a Lie algebra 
homomorphism, we obtain the canonical commutation relation
\begin{equation}
  [ \hat{q}^i, \hat{p}_j ] = \sqrt{-1} \hbar \delta^i_j \,.
  \label{heis-rel}
\end{equation}

This is realized as follows:
Let $\mathcal{H} = L^2( \mathbb{R}^n )$ be the space of $L^2$-functions
on $\mathbb{R}^n$ 
with respect to the Lebesgue measure and define
\begin{align*}
  \hat{q}^i \varphi (q) &= - \sqrt{-1} q^i \varphi (q)\, ,\\
  \hat{p}_j \varphi (q) &= \hbar \frac{\partial \varphi}{\partial q^j} (q)\, ,
\end{align*}
for $\varphi (q) \in L^2( \mathbb{R}^n )$.
Then the relation (\ref{heis-rel}) are satisfied.

\begin{rem}
  From (\ref{heis-rel}),
  $\hat{q}^i$ and $\hat{p}_i$ get commutative when $\hbar \to 0$.
  Such a limit is called the semi-classical limit.
\end{rem}

\begin{df}
  The {\it Heisenberg algebra} $\mathrm{heis}(\mathbb{R}^{2n})$
  is the Lie algebra generated
  by $\hat{q}^1, \dots , \hat{q}^n, \hat{p}_1, \dots , \hat{p}_n$
  (and a center) satisfying the Heisenberg relations.
  In other words, $\mathrm{heis}(\mathbb{R}^{2n})$ is given by
  the central extension
  \[
    0 \longrightarrow \mathbb{R} \longrightarrow
    \mathrm{heis}(\mathbb{R}^{2n}) \longrightarrow \mathbb{R}^{2n}
    \longrightarrow 0\, ,
  \]
  where the extension class is given by $\omega$.
\end{df}

\begin{fact}[Stone, von Neumann (see, for example \cite{M3})]
  $L^2( \mathbb{R}^n )$ is the unique (up to isomorphisms) 
  irreducible representation of $\mathrm{heis}(\mathbb{R}^{2n})$.
\end{fact}

\subsection{Prequantization}

Let $(M,\omega)$ be a symplectic manifold of dimension $2n$
and assume that there exists a prequantum line bundle
$(L, \nabla) \to M$.

As we saw in the previous section, the Lie algebra of the automorphism
group $\mathcal{G}$ of the prequantum bundle $(L,\nabla)$ 
is identified with the space
$C^{\infty}(M)$ of smooth functions on $M$.
Hence we have a representation $\Gamma (L)$ of 
$\mathcal{G}$ and $\mathrm{Lie}\, \mathcal{G} \cong C^{\infty}(M)$.
However this representation is not the desired one.
This is ``too large'' as we see in the following example
(This is why we call this
``prequantization'').

\begin{expl}\label{vecsp}
  We consider the case of symplectic vector space
  $M = T^* \mathbb{R}^n = \mathbb{R}^{2n}$ with the symplectic
  form $\omega = \frac 1{\hbar} \sum_{i=1}^n dq^i \wedge dp_i$,
  as in the previous subsection.
  In this case, the prequantum bundle $L$ is the
  trivial line bundle with connection
  $\nabla = d + \sqrt{-1} \hbar \sum p_i dq^i$.
  Then $\Gamma (L)$ is identified with
  $L^2(M) = L^2(\mathbb{R}^{2n})$.
  Note that
  \[
      \xi_{q^i} = - \hbar \frac{\dr }{\dr p_i}\,, \quad
      \xi_{p_i} = \hbar \frac{\dr}{\dr q^i} \, .
  \]
From this and (\ref{f-hat}), the actions of $q^i$ and $p_i$ are given by
  \begin{align*}
    \hat{q^i} (s) &= - \hbar \frac{\dr }{\dr p_i}s
                     - \sqrt{-1} q^i s \, ,\\
    \hat{p_i} (s) &= \hbar \frac{\dr}{\dr q^i} s
  \end{align*}
  for $s \in L^2(\mathbb{R}^{2n})$.
  This representation is different from the one discussed
  in the previous subsection.
\end{expl}
To get the desired representation, we have to eliminate half of
the parameters $(q^1, \dots , q^n, p_1, \dots , p_n)$.
This is done by choosing a ``polarization''.

\begin{rem}
  We can take $k \omega$ as a symplectic form instead of $\omega$.
  In this case, $L^k$ is a prequantum bundle.
  $\frac 1k$ play the role of the Planck's constant $\hbar$.
  We also call the limit $k \to \infty$ the semi-classical limit.
\end{rem}

\subsection{Polarizations}

\begin{df}
  Extend $\omega$ on $TM \otimes \mathbb{C}$ complex bilinearly.
  A {\it polarization} $P$ is an integrable Lagrangian subbundle
  in $TM \otimes \mathbb{C}$,
  i.e. a complex subbundle $P \subset TM \otimes \mathbb{C}$
  of rank $n$ satisfying $[P,P] \subset P$ and
  $\omega |_P \equiv 0$.
\end{df}

For a polarization $P$, we ``define'' the space of
polarized sections by
\[
  \Gamma_P(L) = \{ s \in \Gamma (L) \, | \, \nabla_{\xi}s=0,\,\,
  \text{for all $\xi \in P$} \} \, .
\]
(This definition is a temporary one.
For some class of polarization, we need to modify the definition.
We will discuss this point later.)

\begin{rem}
The conditions of polarization are necessary for
the integrability condition of $\nabla_{\xi}s=0$.
Since the curvature of $(L,\nabla)$ coincides with $\omega$,
the integrabiliy condition can be written as
\[
  0 = [\nabla_{\xi}, \nabla_{\eta}] s
    = \left( \nabla_{[\xi, \eta]} - \sqrt{-1} \omega(\xi, \eta)
      \right) s
\]
for every $\xi, \eta \in P$.
The integlability of $P$ implies $[\xi, \eta] \in P$ and
the Lagrangian condition implies $\omega(\xi, \eta) =0$.
\end{rem}

There are two important classes of polarizations.

\begin{expl}
Assume that we have a Lagrangian fibration $\pi : (M,\omega) \to B$.
Then the complexified relative tangent bundle
\[
  P := T_{M/B} \otimes \mathbb{C}
     = \ker (d \pi : TM \to TB) \otimes \mathbb{C}
\]
gives a polarization.
This polarization satisfies $\bar{P} = P$.
We call polarizations satisfying this condition
{\it real polarizations}.\footnote{%
In general, Lagrangian fibrations have degenerate fibers.
Thus it is natural to allow real polarizations to be degenerate.
However we ignore singular fibers in this section.}

In this case, the space $\Gamma_P(L)$ of polarized sections
consists of sections of $L$ which are covariantly constant
along fibers.
Note that the restriction of $L$ on each fiber is flat
since each fiber is Lagrangian.
We discuss this case later.
\end{expl}

\begin{expl}
Let $(M,\omega)$ be a K\"ahler manifold.
Then $P := T^{0,1}M \subset TM \otimes \mathbb{C}$ is a polarization.
Such polarization is called a {\it K\"ahler polarization}.
Note that $L$ is holomorphic since the curvature 
$\omega$ is of type $(1,1)$.
K\"ahler polarizations are characterize by
\[
  \begin{cases}
    P \cap \bar{P} = \{0\} \, ,\\
    \omega |_{P \times \bar{P}} > 0 \, .
  \end{cases}
\]
Giving a K\"ahler polarization is equivalent to fixing a compatible
complex structure on $(M,\omega)$.

In this case, $\Gamma_P (L)$ is nothing but the space
$H^0(M,L)$ of holomorphic sections.
\end{expl}

\begin{rem}
The action of $\mathcal{G}$ (or $C^{\infty}(M)$) does not preserve
$\Gamma_P(L) \subset \Gamma (L)$.
In general, the subgroup of $\mathcal{G}$ which preserves
the polarization is very small.
\end{rem}

\begin{expl}
  We consider the case of symplectic vector space again.
  The natural projection
  $\pi : T \mathbb{R}^n \to \mathbb{R}^n$ is a Lagrangian
  fibration.
  Let $P = \ker d \pi$.
  Then $\Gamma_P (L)$ is identified with $L^2(\mathbb{R}^n)$.
  Furthermore, the actions of $q^i$ and 
  $p_i \in C^{\infty}(T \mathbb{R}^n)$
  preserve this subspace.
  In fact, its action is given by
  \[
    \begin{cases}
      \hat{q^i} \varphi (q) = - \sqrt{-1} q^i \varphi (q)\, ,\\
      \hat{p_i} \varphi (q)
        = \hbar \frac{\dr}{\dr q^i} \varphi (q)
    \end{cases}
  \]
  for $\varphi (q) \in L^2(\mathbb{R}^n)$,
  as desired.
\end{expl}

We can take a K\"ahler polarization on 
$T \mathbb{R}^n \cong \mathbb{C}^n$.
In this case, we also have an irreducible representation 
of Heisenberg group on a space of holomorphic functions 
on $\mathbb{C}^n$.
This is called the Bergman-Fock representation.
(See \cite{M3}.)

\begin{expl}[Borel-Weil theory]
  Let $G$ be a compact Lie group, $T$ a maximal torus in $G$,
  $G^{\mathbb C}$ the complexification of $G$,
  and $B \subset G^{\mathbb C}$ a Borel subgroup.
  We denote the flag variety by $X = G/T = G^{\mathbb C} /B$.
  Then every  $G^{\mathbb C}$-equivariant line bundle on $X$ can be 
  determined by a holomorphic character $B \to \mathbb{C}^*$ by
  \[
    L = G^{\mathbb C} \times_B \mathbb{C}.
  \]
  On the other hand, every character of $B$ descend to a character of
  $B/[B,B] \cong T^{\mathbb C}$:
  \[
  \xymatrix{
    B \ar[r] \ar[d] & \mathbb{C}^* \\
    T^{\mathbb C} \ar[ru]& 
  }
  \]
  Hence there is a one to one correspondence
  \[
    \{ \text{$G^{\mathbb C}$-equivariant line bundles on $X$} \}
    \longleftrightarrow
    \{ \text{ holomorphic characters 
    $T^{\mathbb C} \to \mathbb{C}^*$} \}.
  \]

  For each character $\lambda : T^{\mathbb C} \to \mathbb{C}^*$,
  we denote the corresponding line bundle by
  $L_{\lambda} = G^{\mathbb C} \times_B \mathbb{C} \to X$.

\begin{thm}[Borel-Weil]
    If $\lambda$ is a dominant weight, then 
    $H^i(X,L_{\lambda}) = 0$ for $i > 0$ and
    $H^0(X,L_{\lambda})$ is the irreducible representation 
    of $G$ of highest weight $\lambda$.
    Any finite dimensional irreducible representation of $G$
    is given in this way.
  \end{thm}

\end{expl}

Next we discuss the case of real polarizations with compact 
Lagrangian fibers.
Let $\pi : M \to B$ be a Lagrangian torus fibration.
Since $\omega$ is the curvature of $L$,
the restriction $L|_{\pi^{-1} (b)}$ is a flat line bundle
on $\pi^{-1} (b)$ for each $b \in B$.
However $L|_{\pi^{-1} (b)}$ has non-trivial holonomy in general.
Therefore, there is no nontrivial smooth section of $L$ which is 
constant along fibers.
Hence we must change the definition of $\Gamma_P(L)$.

\begin{df}
  A fiber $\pi^{-1}(b)$ is called a {\it Bohr-Sommerfeld fiber} if
  the restriction $(L, \nabla)|_{\pi^{-1}(b)}$ is trivial.
\end{df}

For each $k \in \mathbb{N}$, 
we can take $k \omega$ to be a symplectic form instead of $\omega$.
In this case, $L^k$ is a prequantum bundle.

\begin{df}
  $\pi^{-1}(b)$ is called a {\it Bohr-Sommerfeld fiber of level $k$} 
  (or {\it $k$-Bohr-Sommerfeld fiber}) if
  the restriction $(L^k, \nabla)|_{\pi^{-1}(b)}$ is trivial.
\end{df}

From Proposition \ref{action-angle} and Remark \ref{holonomy},
$\pi^{-1}(b)$ is a Bohr-Sommerfeld fiber of level $k$ if and only
if the action coordinate $(q^1, \dots , q^n)$ of $b$ takes its value
in $\frac 1k \mathbb{Z}^n$.
In particular, Bohr-Sommerfeld fibers appear discretely.

Instead of smooth sections, we consider 
distributional sections of $L$ supported
only on Bohr-Sommerfeld fibers and
are covariantly constant on the support:
\[
  \Gamma_P(L) = \{ s \, | \, \mathrm{supp}\, s \subset
  \text{BS fibers and } \nabla_{\xi}s = 0,\, \xi \in P \}\, .
\]
Since each fiber is connected, covariantly constant sections
on a Bohr-Sommerfeld fiber is unique up to constants.
In particular, the dimension of $\Gamma_P(L)$ is equal to
the number of Bohr-Sommerfeld fibers.

\begin{rem}
There exists a cohomological definition of $\Gamma_P(L)$.
Consider the following complex
\[
  \begin{CD}
    0 @>>> \Gamma (L) @>{d^{\nabla}_P}>> \Gamma (L \otimes P^*)
    @>{d^{\nabla}_P}>> \Gamma (L \otimes \wedge^2 P^*)
    @>{d^{\nabla}_P}>> \dots
  \end{CD}
\]
where $d^{\nabla}_P$ is the restriction of the connection
to $P$.
This is in fact a complex since
\[
  (d^{\nabla}_P)^2 = \omega |_P = 0\, .
\]
Then
\begin{thm}[\'Sniatycki \cite{S}]
\[
  H^i \left(\Gamma (L \otimes \wedge^* P^*),
    d^{\nabla}_P \right) = 0
\]
for $i \ne n$ and
\[
  H^n \left( \Gamma (L \otimes \wedge^* P^*),
  d^{\nabla}_P \right) \cong \Gamma_P (L)
  \quad \text{(``Poincar\'e duality'')}.
\]
\end{thm}
\end{rem}

\subsection{Comparison of the spaces of wave functions for
         K\"ahler and real polarizations}

Assume that we have a Lagrangian fibration
$\pi : (X,\omega) \to B$ of a compact K\"ahler manifold
and a prequantum bundle $L \to X$.
Then we have two spaces of wave functions, i.e. the space
$H^0(X,L)$ of holomorphic sections and
$\Gamma_{T_{X/B}} (L)$.
It is natural to ask whether these spaces are isomorphic or not.
In general, it is difficult to define and compute $\Gamma_{T_{X/B}} (L)$
since there exist degenerate fibers for general Lagrangian fibrations.
However it is shown that these spaces are isomorphic for several cases.

\begin{expl}
Let $X = \mathbb{CP}^1$ with the Fubini-Study metric
and $L = \mathcal{O} (1)$.
The moment map of a natural $S^1$-action is a Lagrangian fibration.
This can be written explicitly as follows.
We identify $\mathbb{CP}^1$ with the unit sphere
\[
  S^2 = \{(x_1,x_2,x_3) \in \mathbb{R}^3 \, |
  \, x_1^2 + x_2^2 + x_3^2 = 1 \,\} \subset \mathbb{R}^3 \, .
\]
Then the moment map is given by
\[
  \pi : S^2 \longrightarrow [-1,1]\, , \quad
  (x_1,x_2,x_3) \longmapsto x_3\, .
\]
It is easy to see that $\pi^{-1}(b)$ satisfies the Bohr-Sommerfeld 
condition of level $k$ if and only if
$b = \frac{2i-k}k$, $i = 0, 1, \dots ,k$
is a lattice point.
Note that 
\[
  \dim H^0(\mathbb{P}^1, \mathcal{O} (k)) = k+1
  = \text{the number of $k$-BS fibers}.
\]
In this case, the isomorphism between two spaces of wave functions
is given by just the restrictions:
\[
  H^0(\mathbb{P}^1, \mathcal{O} (k)) \longrightarrow
  \Gamma_{T_{X/B}}(L^k)\, , \quad
  z_0^i z_1^{k-i} \longmapsto z_0^i z_1^{k-i}|_{\pi^{-1}(\frac{2i-k}k)}\, .
\]
\end{expl}

The same is true for general projective toric manifolds
(with a suitable K\"ahler form).

Guillemin-Sternberg \cite{GS} proved the case of flag varieties.
This is also true for the case of abelian varieties.
We will see this in detail in the next subsection.

\begin{thm}[Andersen \cite{A}] \label{thm-Anderson}
  Let $(X,\omega)$ be a compact K\"ahler manifold, $L \to X$ a 
  prequantum bundle.
  Assume that we have a Lagrangian fibration $\pi : X \to B$
  with no degenerate fiber.
  Then
  \[
    \dim H^0(X,L^k) = \dim \Gamma_{T_{X/B}} (L^k)
  \]
  for large $k$.
\end{thm}

\begin{rem}
  In general, the restriction of holomorphic sections to BS fibers
  does not give covariantly constant sections.
\end{rem}

\begin{rem}
Very few compact K\"ahler manifolds admit Lagrangian fibration
without degenerate fibers.
However it is expected that this theorem holds for more general cases
such as K3 surfaces.
We discuss this point in section 4 again.
\end{rem}

The dimension of $H^0(X,L^k)$ is given by Riemann-Roch formula
and vanishing theorem of cohomologies.
Thus what we need to do is to count the number of Bohr-Sommerfeld fibers.
The idea is to use its ``mirror'', i.e. a dual torus
fibration.

\begin{proof}
From Riemann-Roch theorem,
\[
  \dim H^0(X,L^k) = \int_X ch (L^k) \hat{A}(TX)
                  = \int_X \exp (k \omega) \hat{A}(TX)\, .
\]
Now we have
\[
  0 \longrightarrow T_{X/B} \longrightarrow TX \longrightarrow
  \pi^* TB \longrightarrow 0.
\]
By using the symplectic form, $T_{X/B}$ can be identified with
$\pi^* T^*B$.
Since $T^*B$ and $TB$ contains lattices of maximal rank,
these are flat.
Then
\[
  \hat{A}(TX) = \hat{A}(T_{X/B}) \hat{A}(\pi^* TB)
              = \pi^* \bigl( \hat{A} (T^*B) \hat{A} (TB) \bigr) = 1\, . 
\]
Consequently we have
\[
  \dim H^0(X,L^k) = \int_X k^n \omega^n/ n!
                  = k^n \cdot \mathrm{vol}(X,\omega)\, .
\]

Next we calculate the number of Bohr-Sommerfeld fibers.
First we define a dual torus fibration of $\pi :X \to B$.
Let $T^*B/\Lambda \to B$ be the Lagrangian fibration associated to
$\pi$, as in the previous section.
We denote the dual lattice of $\Lambda$ by
$\Lambda^* = \mathrm{Hom}(\Lambda, \mathbb{Z}) \subset TB$.
Then $\check{\pi} : TB/\Lambda^* \to B$ is a dual torus fibration.

Recall that the dual torus $\check{T}^n$ of a torus $T^n$ parameterizes
flat line bundles $E \to T^n$.
In our situation, for each $b \in B$, we can associate a flat line bundle
$L^k|_{\pi^{-1}(b)}$ on $\pi^{-1}(b)$.
This defines a section $\lambda_k : B \to TB/\Lambda^*$ 
of the dual torus fibration.
On the other hand, the zero section $\lambda_0$ corresponds to the
trivial bundle on $X$.
By definition, $\pi^{-1}(b)$ is a Bohr-Sommerfeld fibers of level $k$
if and only if $\lambda_k(b) = \lambda_0(b)$.
This implies 
\[
  \dim \Gamma_{T_{X/B}} (L^k) = \#( \lambda_0(B) \cap \lambda_k(B))\, .
\]
  
$\lambda_k$ can be written explicitly as follows.
Let $(x^1, \dots , x^n, y_1, \dots , y_n)$ be the action-angle coordinate
constructed in section 2 and take the dual coordinate
$(x^1, \dots , x^n, v^1, \dots , v^n)$.
Then $\lambda_k$ is given by
\begin{equation}
  \lambda_k (x^1, \dots , x^n) = (x^1, \dots , x^n, kx^1, \dots , kx^n)\, .
  \label{dual-section}
\end{equation}
From this expression, we see that $\lambda_k$ and $\lambda_0$
intersect transversely and positively (under a suitable orientation).
Therefore the number of Bohr-Sommerfeld fibers of level $k$
coincides with the intersection number of
$\lambda_k(B)$ and $\lambda_0(B)$.
Put $\alpha = dv^1 \wedge \dots \wedge dv^n \in \Omega^n(TB/P)$.
This is the Poincar\'e dual of $\lambda_0(B)$.
\begin{align*}
  \dim \Gamma_{T_{X/B}} (L^k)
  &= \int_{\lambda_k(B) \cap \lambda_0(B)}1\\
  &= \int_{\lambda_k(B)} \alpha
  = \int_B \lambda_k^* \alpha \, .
\end{align*}
From (\ref{dual-section}), we have
\[
  \lambda_k^* \alpha = k^n dx^1 \wedge \dots \wedge dx^n
  = k^n \pi_* \omega^n/n!,
\]
hence
\begin{align*}
  \dim \Gamma_{T_{X/B}} (L^k) &= k^n \int_B \pi_* \omega^n/n!\\
  &= k^n \int_X \omega^n/n!
  = k^n \cdot \mathrm{vol}(X,\omega) \, .
\end{align*}
\end{proof}

\subsection{Geometric quantization for abelian varieties}

Let $X= \mathbb C^n /( \Omega \mathbb Z^n + \mathbb Z^n)$ be 
an abelian variety of complex dimension $n$,
where $\Omega$ is an $n \times n$ complex symmetric matrix whose 
imaginary part $\mathrm{Im}\,\Omega$ is positive definite,
and $L$ an ample line bundle on $X$.
We assume that $L$ is symmetric of degree 1, for simplicity.
Recall that $L$ is said to be symmetric if $(-1_X)^*L \cong L$, here 
$-1_X:X \to X$ is the inverse morphism of $X$.
Explicitly, $L$ is given by
\begin{equation}
  L = ( \mathbb C^n \times \mathbb C )/(\Omega \mathbb Z^n + \mathbb Z^n)\, ,
  \label{dfL}
\end{equation}
where the action of $\Omega \mathbb Z^n + \mathbb Z^n$ 
on $\mathbb C^n \times \mathbb C$ is
\[
  (z,\zeta) \mapsto (z+\lambda, 
           e^{\pi {}^t \bar{\lambda} (\mathrm{Im}\,\Omega)^{-1}z
               + \frac{\pi}2 {}^t \bar{\lambda} (\mathrm{Im}\,\Omega)^{-1} \lambda}
            \zeta)
\]
for $\lambda \in \Omega \mathbb Z^n + \mathbb Z^n$.
Any ample bundle of degree 1 is a pull back of $L$ by some translation
of $X$.
From the definition of $L$, 
\[
  h_0 = \exp(-\pi {}^t z (\mathrm{Im}\,\Omega)^{-1} \bar z)
\]
defines a Hermitian metric on $L$ whose first Chern form
\[
  c_1(L,h_0) = - \frac{\sqrt{-1}}{2 \pi} \partial \bar\partial \log h_0
             =  \frac{\sqrt{-1}}2 \sum_{i,j} h_{ij} dz^i \wedge d \bar z^j
             =: \omega_0
\]
gives a flat K\"ahler metric on $X$, 
here we put $(h_{ij}) = (\mathrm{Im}\,\Omega)^{-1}$.

Take the following two Lagrangian tori
\[
  X^+ = \{\Omega x \, | \, x \in \mathbb{R}^n/ \mathbb{Z}^n\}\, ,\,\, 
  X^- = \{ y \, | \, y \in \mathbb{R}^n/ \mathbb{Z}^n\} \subset X
\]
and consider the natural projection 
\begin{equation}
  \pi : X = X^+ \times X^- \longrightarrow X^- \, , \quad
  z = \Omega x + y \longmapsto y \, .
  \label{fib}
\end{equation}

Then we have two polarizations and hence two
spaces of wave functions.

Next we define finite Heisenberg groups in a similar way as in the case of
vector space $\mathbb{R}^{2n}$.
For each $k \in \mathbb N$, let $X_k$ (resp. $X^{\pm}_k$) 
be the subgroup of points 
of $X$ (resp. $X^{\pm}$) of order $k$.
Then $X^{\pm}_k \cong \frac 1k \mathbb Z^n / \mathbb Z^n$ 
and $X_k = X^+_k \times X^-_k$.  
It is easy to see that $X_k$ can be characterized by
\[
  X_k = \{ w \in X |\, \tau^*_w L^k \cong L^k \}\, ,
\]
where $\tau_w : X \to X$ will denote the translation $\tau_w (z) = z + w$.
Then there exists a central extension $\mathcal{G} (L^k) \subset 
\mathrm{Aut}(X,L^k)$ of $X_k$:
\[
  \begin{CD}
    1 @>>> \mathbb C^* @>>> \mathcal{G} (L^k) @>>> X_k @>>> 0\, ,
  \end{CD}
\]
where $\mathbb C^*$ acts on $L^k$ by multiplications on each fiber.
The {\it finite Heisenberg group} $G_k$ is defined by
\[
  \begin{CD}
    1 @>>> \mathbb C^* @>>> \mathcal{G} (L^k) @>>> X_k @>>> 0\\
      @.   \cup     @.    \cup             @.   \Vert @.\\
    1 @>>> \mu_k   @>>>   G_k      @>>> X_k @>>> 0
  \end{CD}
\]
where $\mu_k$ is the group of $k$-th roots of 1.
The group law is 
\[
  (c_1,a_1,b_1) \cdot (c_2,a_2,b_2) = 
       (e^{2 \pi \sqrt{-1} \xt b_2 a_1}c_1 c_2, a_1+a_2, b_1+b_2)
\]
for $(c_1,a_1,b_1),(c_2,a_2,b_2) \in G_k \cong \mu_k \times \frac 1k \mathbb Z^n
\times \frac 1k \mathbb Z^n$.


\subsubsection*{K\"ahler polarization}

$G_k$ acts naturally  on $H^0(X,L^k)$.
Note that the Hermitian metric $h_0$ is invariant under the $G_k$-action.
This means that the induced $G_k$-action on $H^0(X,L^k)$ is unitary.

\begin{thm}
  $H^0(X,L^k)$ is a unique (up to isomorphism) 
  irreducible representation of $G_k$, called the Heisenberg representation.
\end{thm} 

Holomorphic sections of $L^k$ can be explicitly written 
by using theta functions.
We take a basis of $H^0(X,L^k)$ defined by
\begin{equation}
  s_{b_i} = s_i = 
   C_{\Omega} k^{- \frac n4} \cdot \theta_0(\Omega, z)^k \cdot \tht{-b_i}{z}\, ,
  \label{theta}
\end{equation}
where we denote $X_k^- = \{ b_i\}_{i = 1, \dots , k^n}$ and 
\begin{align*}
  C_{\Omega} &= 2^{\frac n4}(\det (\mathrm{Im}\,\Omega))^{\frac 14}\, ,\\
  \theta_0(\Omega, z) &=  
     \exp \left(\frac {\pi}2 \xt z (\mathrm{Im}\,\Omega)^{-1}z \right)\, ,\\ 
  \vartheta\genfrac[]{0pt}{}{a}{b}(\Omega,z) &=  
     \sum_{l \in \mathbb Z^n} 
       e \left( \frac 12 \xt (l+a)\Omega (l+a) + \xt (l+a)(z+b) \right)\, ,\\
  e(t) &=  \exp (2 \pi \sqrt{-1}t)\, .
\end{align*}
By definition, the action of $G_k$ on $H^0(X,L^k)$ is given by
\[
  s(z) \longmapsto c \cdot \tau^*_{\Omega a + b} s 
                    =  c \cdot s(z + \Omega a + b)
\]
for $(c,a,b) \in G_k$.
Put $\{a_j\} = X^+_k$.
By direct computation, we have
\begin{equation}
  \begin{split}
    \rho_k(1,a_j,0)s_{b_i} &= e^{2 \pi \sqrt{-1} \xt a_j b_i} s_{b_i} \, ,\\
    \rho_k(1,0,b_j)s_{b_i} &= s_{b_i - b_j}\, .
  \end{split}
  \label{rep}
\end{equation}

From (\ref{rep}), we can show that $H^0(X,L^k)$ is an 
irreducible representation of $G_k$.


\subsubsection*{Real polarizations}

Next we consider the real polarization defined by (\ref{fib}).
In this case, Bohr-Somerfeld fibers are characterized as follows:
\begin{prop}
  For $b \in X^-$, $\pi^{-1}(b) = X^+ \times \{b\}$ satisfies the Bohr-Sommerfeld
  condition of level $k$ if and only if $b \in X^-_k$.
\end{prop}

In fact, for each $b_i \in X^-_k$,
\[
  \sigma_i (x) = \exp \left( \frac{k \pi}2
        \xt (\Omega x + b_i) (\mathrm{Im}\,\Omega)^{-1} (\Omega x + b_i)
            - \sqrt{-1} k \pi  \xt x \Omega x \right)
\]
defines a covariantly constant section of $L^k|_{\pi^{-1}(b_i)}$.

In particular the number of Bohr-Somerfeld fibers of level $k$ coincides 
with  $\dim \, H^0(X,L^k) = k^n$.

The basis $s_i \in H^0(X,L^k)$ can be reconstructed from $\sigma_i$'s by
using the Bergman kernel.
The Bergman kernel $\Pi_k(z,w)$ is the integral kernel
of the orthogonal projection
\[
  L^2(X,L^k) \longrightarrow H^0(X,L^k)\, .
\] 

\begin{prop} \label{Berg}
  Let $\Pi_k(z,w)$ be the Bergman kernel of the 
  Hermitian line bundle $(L^k,h_0)$.
  Then
  \[
    \int_{X^+ \times \{b_i\}} \Pi_k(z,x) \sigma_i(x) \, dx
      = C_{\Omega}'k^{\frac n4} s_i(z)\, ,
  \]
  where
 \[
    C_{\Omega}'
      = \frac {2^{\frac n4} \sqrt{-1}^{\frac n2}\det (\mathrm{Im}\, \Omega)^{\frac n4}}
              {(\det \bar{\Omega})^{\frac n2}}\, .
  \]
\end{prop}

\begin{proof}
We will show that $s_1,\dots ,s_{k^n}$ define an orthonormal
basis of $H^0(X,L^k)$ in the next subsection.
Then the Bergman kernel is given by
\[
  \Pi_k(z,w) = \sum_{i=1}^{k^n} s_i(z) s_i(w)^*.
\]
From this expression, it suffices to show that
\begin{equation}
  \int_{X^+ \times \{b_i\}} (\sigma_i , s_j)_{h_0} dx
         = C_{\Omega}' k^{\frac n4} \delta_{ij}\, .
  \label{inner-product}
\end{equation}
By direct computation, we have
\[
  (\sigma_i , s_j)_{h_0} =
    C_{\Omega} k^{\frac n4} \sum_{l \in {\mathbb Z}^n}
    e \left( - \frac k2 \left( \lt \left( x + \frac lk \right) \bar{\Omega}
    \left( x + \frac lk \right) \right) - \xt l (b_i -b_j) \right),
\]
and we obtain (\ref{inner-product}) by integrating this.
\end{proof}



\newpage
\section{Special Lagrangian fibrations and mirror symmetry}

Mirror symmetry is a duality in string theory.
Mathematically, this is regarded as a duality between 
symplectic geometry on a Calabi-Yau manifold $M$ and 
complex geometry on another Calabi-Yau manifold $W$.
Homological mirror symmetry conjectured by Kontsevich \cite{Ko}
claims the equivalence of the Fukaya category on $M$ and
the derived category of coherent sheaves on $W$.
Strominger-Yau-Zaslow \cite{SYZ} proposed mirror symetry via
special Lagrangian fibrations.
This picture gives a geometric construction of mirror manifolds
and correspondence of two categories.
This works successfully for abelian varieties and
theta functions play an important role again.

\subsection{Special Lagrangian submanifolds}

\begin{df}
  A {\it Calabi-Yau manifold} $X$ is a K\"ahler manifold
  with trivial canonical line bundle: $K_X \cong \mathcal{O}_X$.
\end{df}

An important fact is:
\begin{thm}[Yau \cite{Y}]
 If $X$ is a compact K\"ahler manifold with
 $c_1(X) = c_1(K_X^{-1}) = 0$, then
 each K\"ahler class contains a unique Ricci-flat K\"ahler metric.
\end{thm}

Let $X$ be a compact Calabi-Yau manifold of complex dimension $n$.
By definition,
$X$ carries a non-vanishing holomorphic $n$-form
$\Omega \in H^0(X, K_X)$.
Recall that the Ricci form of a K\"ahler form $\omega$ is given by
$\mathrm{Ric}(\omega) = - \dd \log \omega^n$.
Hence $\omega$ is Ricci-flat if and only if
$\omega^n = c \Omega \wedge \bar{\Omega}$ for some constant $c$.

Now we fix a Ricci-flat metric $\omega$ and normalize $\Omega$ so that
\[
  \frac{\omega^n}{n!} = (-1)^{\frac{n(n-1)}2} 
  \left( \frac{\sqrt{-1}}2 \right)^{\frac n2} 
  \Omega \wedge \bar{\Omega}\, .
\]
This condition is satisfied for 
$\omega = \frac{\sqrt{-1}}2 \sum_{i=1}^n dz^i \wedge d \bar{z}^i$
and $\Omega = dz^1 \wedge \dots \wedge dz^n$.
Such $\Omega$ is unique up to $S^1$.

\begin{df}
 Let $S \subset X$ be an $n$-dimensional submanifold.
 $S$ is called a {\it special Lagrangian submanifold} if
 \[
   \omega|_S = \im (e^{\sqrt{-1}\theta} \Omega)|_S = 0
 \]
 for some $\theta \in \mathbb{R}$.
\end{df}

An important property of special Lagrangian submanifolds is
the following:
\begin{thm}[Harvey-Lawson \cite{HL}]
  Every special Lagrangian submanifold is volume minimizing in its
  homology class.
\end{thm}

Note that the special Lagrangian condition is a first order 
differential equation.
On the other hand, the condition of minimal submanifolds
is second order.
In addition, minimal submanifolds are not neccesarily volume
minimizing.
Such a situation is similar to that of ASD connections and
Yang-Mills connections in gauge theory.

Recall that the normal bundle of a Lagrangian submanifold $S$
is identified with the cotangent bundle $T^*S$ by using 
the symplectic structure.
Under this identification, every small deformations of $S$ is given as
a graph of a 1-form on $S$.
Then the graph is also Lagrangian if and only if 
the 1-forms is closed.

\begin{thm}[McLean \cite{Mc}]
 Under the above identification, infinitesimal deformations of
 a special Lagrangian submanifold are given by harmonic 1-forms:
 \[
   d \alpha = d^* \alpha = 0\, .
 \]
 Furthermore, every infinitesimal deformation is unobstructed.
 In particular, the tangent space of special Lagrangian submanifolds
 is identified with $H^1(S,\mathbb{R})$.
\end{thm}

\begin{df}
  A fibration $\pi : X \to B$ is called a special Lagrangian fibration 
  if each smooth fiber is special Lagrangian. 
\end{df}

Since each fiber $\pi^{-1}(b)$ is a Lagrangian torus,
\[
  T_b B \overset{\omega}{\cong} 
  T^*_x \pi^{-1}(b) \cong H^1(\pi^{-1}(b), \mathbb{R})
\]
for $x \in \pi^{-1}(b)$.
This implies that $B$ can be considered as 
(a component of) the moduli space of 
special Lagrangian tori in $X$.

\begin{expl}
  Let $X = \mathbb{C}^n/(\tau \mathbb{Z}^n +  \mathbb{Z}^n)$
  be an abelian variety with a flat metric
  $\omega = \frac{\sqrt{-1}}2 \sum h_{ij} dz^i \wedge d \bar{z}^j$,
  where $\tau$ is a $n \times n$ matrix with positive definite imaginary
  part and $(h_{ij}) = (\im \tau)^{-1}$.
  Then
  \[
    \pi : X \longrightarrow T^n , \quad 
    \tau x + y \longmapsto x
  \]
  is a special Lagrangian fibration.
\end{expl}

\begin{expl}
  Let $X$ be a K3 surface and assume that $\pi : X \to \cp{1}$ is
  an elliptic fibration with respect to a complex structure $I$.
  We fix a Ricci-flat K\"ahler metric $g$.
  Since $X$ is hyperK\"ahler, there exist compatible complex structures
  $J$ and $K$ satisfying 
  $I^2= J^2 = K^2 = -1$ and $IJ=K$.
  We denote the corresponding K\"ahler form by $\omega_I$, $\omega_J$
  and $\omega_K$ respectively.
  Then
  \begin{align*}
    \Omega_I = \omega_K + \sqrt{-1} \omega_J &\, 
    \text{ is holomorphic with respect to $I$,}\\
    \Omega_J = \omega_I + \sqrt{-1} \omega_K &\, 
    \text{ is holomorphic with respect to $J$,}\\
    \Omega_K = \omega_J + \sqrt{-1} \omega_I &\, 
    \text{ is holomorphic with respect to $K$.}
  \end{align*}
  Since $\pi$ is holomprphic with respect to $I$,
  we have $\Omega_I|_{\mathrm{fiber}} = 0$.
  This implies that $\pi : X \to S^2$ is a special Lagrangian
  fibration with respect to $J$.
  
  In this case, we know the structure of singular fibers
  from the result for elliptic surfaces
  by Kodaira \cite{Kod}.
\end{expl}

\subsection{A local model of mirror symmetry}

Homological mirror symmetry conjecture 
states the equivalence of the Fukaya category of $M$ and  
the derived category of coherent sheaves on $W$.
Roughly, this means that
each pair $(S,L)$ of a special Lagrangian submanifold in $M$
and a flat line bundle on $S$ corresponds to a coherent sheaf on $W$.
In this subsection, we see this correspondence from the point of view of 
special Lagrangian fibrations,
following Leung-Yau-Zaslow \cite{LYZ} and Leung \cite{L}
(see also Hitchin \cite{Hi}).

\begin{conj}[Strominger-Yau-Zaslow \cite{SYZ}]
If $(M,W)$ is a mirror pair,
then there exist special Lagrangian  torus fibrations
\begin{align*}
 \pi & : M \longrightarrow  B\, ,\\
 \check{\pi} &: W  \longrightarrow B
\end{align*}
satisfying the following condition:
there exists an open dense subset $B_0 \subset B$ such that
$\pi^{-1}(b)$ and $\check{\pi}^{-1}(b)$ are smooth and 
dual to each other
for each $b \in B_0$.
\end{conj}

This picture gives a concrete description of
the correspondence of (special) Lagrangian submanifolds in $M$ 
and coherent sheaves on $W$.
In this subsection, we see the correspondence of  
holomrphic condition and special Lagrangian condition 
in a local model.

Let $B$ be a small ball in $\mathbb{R}^n$ with coordinate
$(x^1, \dots , x^n)$.
Then every Lagrangian torus fibration has the form
\[
  \pi : M = T^*B/\Lambda \longrightarrow B
\]
for some $\mathbb{Z}^n$-bundle $\Lambda \subset T^*B$.
We denote the standard symplectic form by
$\omega = \sum dx^i \wedge dy_i$.
As in the previous section, we construct the dual torus fibration by
\[
  \check{\pi} : W = TB/\Lambda^* \longrightarrow B \, ,
\]
where $\Lambda^* = \mathrm{Hom}(\Lambda, \mathbb{Z})$ is the
dual lattice bundle of $\Lambda$.
Then $W$ has a complex structure induced from
\[
  TB \longrightarrow \mathbb{C}^n\, , \quad
  \sum_{i=1}^n y^i \frac{\dr}{\dr x^i} \longmapsto 
  (z^1 = x^1 + \sqrt{-1} y^1, \dots)\, ,
\]
where $y^i$ is the dual coordinate of $x^i$.
Note that $\Omega = dz^1 \wedge \dots \wedge dz^n$ is a 
holomorphic $n$-form on $X$.

\begin{rem}
  Here we assume that the B-field is zero.
  If we consider non-zero B-fields, we obtain more general 
  complex structures.
\end{rem}

Next we introduce K\"ahler metrics on $M$ and $W$.
From the definition of the complex structure on $W$,
every K\"ahler metric and its K\"aher form have the form
\begin{align*}
  \check{g} &= \sum g_{ij}(dx^i dx^j + dy^i dy^j)\, ,\\
  \check{\omega} &= \frac{\sqrt{-1}}2 \sum g_{ij} dz^i \wedge d \bar{z}^j\, .
\end{align*}
The induced metric on $M$ is written as
\[
  g = \sum (g_{ij} dx^i dx^j + g^{ij} dy_i dy_j)\, ,
\]
where $(g^{ij}) = (g_{ij})^{-1}$.
Then we have a compatible almost complex structure $J$
on $M$ defined by
\[
  \omega( \xi, \eta) = g (J \xi , \eta)\, .
\]

Recall that each fiber of $\pi$ and $\check{\pi}$ have a natural 
affine structure.
Hereafter we assume that $g$ is constant along fibers with respect to 
the natural affine structure (semi-flat condition).
Then the K\"ahler condition ($d \omega = 0$) implies that
\[
  g_{ij} = \frac{\dr^2 \phi}{\dr x^i \dr x^j}
\]
for some $\phi = \phi (x) \in C^{\infty}(B)$.
Since 
$(g_{ij}) = \left(\frac{\dr^2 \phi}{\dr x^i \dr x^j} \right)$
is non-degenerate, we can take a new coordinate
$(x_1, \dots , x_n)$ of $B$ satisfying
\[
  x_i = \frac{\dr \phi}{\dr x^i}\, , \quad i= 1, \dots , n.
\]
Remark that the Jacobian is given by
$\frac{\dr x_i}{\dr x^j} = g_{ij}$.
We denote the Legendre transform of $\phi$ by
\[
  \psi = \sum_{i=1}^n x^i x_i - \phi \, .
\]
Then 
\[
  x^i = \frac{\dr \psi}{\dr x_i} \quad 
  \text{and} \quad
  g^{ij} = \frac{\dr^2 \psi}{\dr x_i \dr x_j}\, .
\]
Hence we have
\begin{align*}
  g &= \sum g^{ij} (dx_i dx_j + dy_i dy_j),\\
  \omega &= \frac{\sqrt{-1}}2 \sum g^{ij} dz_i \wedge d \bar{z}_j \, ,
\end{align*}
where $z_i = x_i + \sqrt{-1} y_i$.
In particular $J$ is integrable.

The Ricci-flat condition for $g$ is equivalent to
the Monge-Amp\`ere equation
\[
  \det \left( \frac{\dr^2 \phi}{\dr x^i \dr x^j} \right) = C
\]
for some constant $C$.
Since $\det (g^{ij}) = \det (g_{ij})^{-1}$,
$g$ is Ricci-flat if and only if $\check{g}$ is Ricci-flat.

Namely, in this setting, mirror of $M$ is nothing but its Legendre
transform.

We summarize the above discussion.
$M = T^*B/ \Lambda$ has a standard symplectic structure
while $W = TB/\Lambda^*$ has a natural complex structure.
Each compatible almost complex structure $J$
(or, equivalently, a metric $g$)
on $M$ determines a metric $\check{g}$ (or $\check{\omega}$) on $W$.
Under the semi-flat condition,
closedness of $\check{\omega}$ corresponds to the integrability condition
of $J$.
Moreover, Ricci-flat condition for $g$ is equivalent to that for 
$\check{g}$.

\begin{rem}
  The semi-flat condition is quite strong.
  In fact, few compact Calabi-Yau manifolds admit semi-flat
  Ricci-flat metrics.
  Nevertheless this local model is still important
  because it is expected that Calabi-Yau manifolds 
  posses structures close to this semi-flat model
  near the large complex structure limit 
  (i.e when it is close to the most degenerate Calabi-Yau's).
  Gross-Wilson \cite{GW} analyzed the case of K3 surfaces. 
\end{rem}

Next we consider the correspondence of special 
Lagrangian submanifolds in $M$ and 
coherent sheaves on $W$ .

\begin{df}[Leung-Yau-Zaslow \cite{LYZ}]
\begin{enumerate}
  \item A {\it supersymmetric A-cycle} is a pair $(S,L)$
    of a special Lagrangian submanifold $S$ and a flat
    line bundle on it.
  \item A {\it supersymmetric B-cycle} is a pair $(C,E)$
    of a complex submanifold $C$ and a holomorphic vector
    bundle $E \to C$ 
    with a unitary connection $A$ satisfying 
    $\im e^{\sqrt{-1}\theta} (\omega + F_A)^m = 0$
    for some $\theta \in \mathbb{R}$,
    where $F_A$ is the curvature of $E$ and
    $m = \dim_{\mathbb{C}}C$.
    We call holomorphic vector bundles satisfying this condition
    {\it Generalized Hermitian Yang-Mills}.
\end{enumerate}
\end{df}

\begin{rem}
  A holomorphic vector bundle $(E,A) \to C$ is said to be
  {\it Hermitian-Yang-Mills} if 
  $ \Lambda F_A = \lambda Id$ for some constant $\lambda$,
  where $\Lambda$ is the adjoint operator of $\omega \wedge$.
  In this case, $F_A$ satisfies the generalized 
  Hermitian-Yang-Mills condition.
\end{rem}

First we consider the case that $S \subset M$ is a
special Lagrangian section of $\pi : M \to B$.
$S$ can be given as a graph of 1-form 
$\alpha = \sum A_i dx^i$ on $B$.
The Lagrangian condition is equivalent to 
$d \alpha = 0$.
We set
\begin{equation}
  \frac{\dr A_i}{\dr x^j} = \frac{\dr A_j}{\dr x^i} = F_{ij}\, .
  \label{closed}
\end{equation}
Since $z_i = x_i + \sqrt{-1} y_i = 
\frac{\dr \phi}{\dr x^i} + \sqrt{-1} y_i$,
we have
\[
  \alpha^* dz_i = \frac{\dr}{\dr x_j} 
  \left( \frac{\dr \phi}{\dr x^i} + \sqrt{-1} A_i \right) dx^j
  = (g_{ij} + \sqrt{-1} F_{ij})dx^j\, .
\]
Hence
\[
  \alpha^* \Omega = \det (g_{ij} + \sqrt{-1} F_{ij})
  dx^1 \wedge \dots \wedge dx^n.
\]
The special Lagrangian condition 
$\im e^{\sqrt{-1}\theta} \Omega = 0$ becomes
\begin{equation}
  \im e^{\sqrt{-1}\theta} \det (g_{ij} + \sqrt{-1} F_{ij}) = 0\, .
  \label{special}
\end{equation}
Then the corresponding supersymmetric B-cycle is given by
$C = W$ and $E$ is a line trivial bundle with the connection
$\nabla_A = d + \sqrt{-1} \sum A_i(x)dy^i$.
Note that $(E,A)$ is holomorphic if and only if
\[
  F_A^{0,2} = \frac 14 \sum F_{ij} d \bar{z}^i \wedge d \bar{z}^j
  = 0\, .
\]
This follows from (\ref{closed}) (i.e. the closedness of $\alpha$).
Furthermore the generalized Hermitian-Yang-Mills condition is
equivalent to the special Lagrangian condition (\ref{special}).

Next we consider the case that $S$ is a multisection, i.e.
the restriction $\pi |_S : S \to B$ is a $r$-fold covering.
In this case, the corresponding supersmmetric B-cycle is a
holomorphic vector bundle $E \to W$ of rank $r$.
To see this, we consider the following $r$-fold covering:
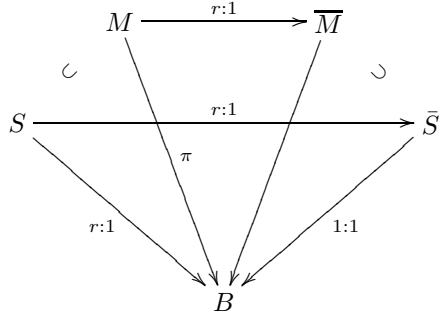
\[
  \xymatrix{
    & M \ar[rr]^{r:1} \ar[dddr]^{\pi} & 
   & \overline{M} \ar[dddl] & \\
   S \ar[ddrr]_{r:1} \ar[rrrr]^{r:1} 
   \ar@{}[ur]|{\text{\rotatebox[origin=c]{45}{$\subset$}}}
   & & & & \bar{S}\ar[ddll]^{1:1} 
   \ar@{}[ul]|{\text{\rotatebox[origin=c]{-45}{$\supset$}}}\\ 
   & & & & \\
    & & B & &
   }
\]
Then the problem is reduced to the case of single sections.
$\bar{S}$ corresponds to a line bundle $L \to \widetilde{W}$.
$E$ is defined by
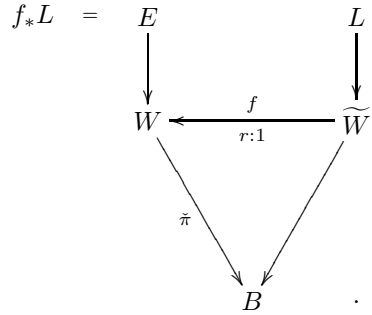
\[
   \xymatrix{
   f_*L \ar@{}[r]|= & E \ar[d] & & L \ar[d]\\
   & W \ar[ddr]_{\check{\pi}} & &
   \widetilde{W} \ar[ll]^{r:1}_f \ar[ddl] \\
   & & &  \\
    & & B & .
   }
\]
Since $ f : \widetilde{W} \to W$ is a $r$-fold covering,
$E$ is a vector bundle of rank $r$.
\subsection{Mirror symmetry for abelian varieties}

The correspondence discussed in the previous subsection
works in the case of abelian varieties.
For elliptic curves,
Polishchuk-Zaslow \cite{PZ} proved an equivalence of 
the derived category
$D^b(W)$ of coherent sheaves on $W$ and a modified version 
$\mathcal{F}_0(M)$ of Fukaya category
from this view point.
For higher dimensional case,
Fukaya \cite{F} constructed a functor from the $A_{\infty}$ category
of affine Lagrangian submanifolds in $M$
to $D^b(W)$ for abelian varieties.
We recall these theory very briefly.

Let 
$M = \mathbb{R}^{2n}/\mathbb{Z}^{2n}$ be a $2n$ diemsional torus
with a constant complexified symplectic form
$\omega^{\mathbb{C}} =\omega + \sqrt{-1}\beta$, where
\begin{align*}
  \omega &= \sum dx^i \wedge d y^i\, ,\\
  \beta &= \sum b_{ij} dx^i \wedge d y^j\, ,
\end{align*}
and consider the following Lagrangian fibration
\[
  \pi : M \longrightarrow B=T^n, \quad 
  (x,y) \longmapsto x\, .
\]

Set $\tau = (b_{ij} + \sqrt{-1} \delta_{ij})$.
Then the mirror is given by
\[
  W = \mathbb{C}^n / \tau \mathbb{Z}^n + \mathbb{Z}^n
\]
with
dual torus fibration $\check{\pi} : W \to B$.

\begin{thm}[Fukaya \cite{F}, Polischuk \cite{P}]
For each pair $(S,L)$ of an 
affine Lagrangian submanifold $S$ in $M$ 
and a flat line bundle $L \to S$, we can associate an object
$E=E(S,L)$ of $D^b(W)$.
\end{thm}

The construction in the previous subsection can be rephrased
by using ``the Poincar\'e bundle''.
Let $p : M \times_B W \to W$ be the natural projection.
Then the Poincar\'e bundle is a line bundle
\[
  \mathcal{P} \to M \times_B W
\]
such that, for $\xi \in W$, 
\[
  \mathcal{P}|_{p^{-1}(\xi)} \to \pi^{-1}( \check{\pi}(\xi))
\]
is the flat line bundle corresponding to $\xi
\in \check{\pi}^{-1}( \check{\pi}(\xi))
= (\pi^{-1}( \check{\pi}(\xi)))^{\vee}$.

Let $S \subset M$ be an affine Lagrangian submanifold and
$L \to S$ a line bundle.
Assume that $\pi|_S : S \to B$ is an unramified covering.
Let
\begin{align*}
  p_S &: S \times_B W \to S\, ,\\
  p_W &: S \times_B W \to W\, ,\\
  (i \times id) &: S \times_B W \to M \times_B W
\end{align*}
be natural maps.
Then $E(S,L)$ is defined by
\[
  E(S,L) = (p_W)_* \bigl( (i \times id)^* \mathcal{P} \otimes
            p_S^* L \bigr)
\]
(``Fourier-Mukai transform'' of $(S,L)$).

\begin{prop}[Polischuk \cite{P}]
  $E(S,L)$ is holomorphic if and only if
  the curvature of $L$ coincides with
  $\omega^{\mathbb{C}}|_S =0$.
\end{prop}

Newt we see the correspondence of morphisms.
Morphisms in Fukaya category are given by Floer homologies
$HF\bigl( (S_1,L_1), (S_2,L_2) \bigr)$.
Roughly speaking, Flore homologies are generated by intersection points
of Lagrangian submanifolds.
Here we are interested in the case that Lagrangian submanifolds
intersect transversally.

\begin{fact}[Fukaya]
If $S_1$ and $S_2$ are transverse,
\[
  HF \bigl( (S_1,L_1), (S_2,L_2) \bigr)
  \cong \mathbb{C}^{\# S_1 \cap S_2} 
  \otimes \mathrm{Hom}(L_1,L_2)\, ,
\]
where $\mathrm{Hom}(L_1,L_2)$ is homomorphisms of vector spaces
underlying the local systems at $S_1 \cap S_2$.
\end{fact}

\begin{thm}[Fukaya \cite{F}]
Let $S_1$ and $S_2$ be affine Lagrangian submanifolds in $M$
(not necessarily transverse) and $L_1$, $L_2$ flat line bundles
on $S_1$ and $S_2$ respectively.
Then
\[
  HF \bigl( (S_1,L_1), (S_2,L_2) \bigr)
  \cong \mathrm{Ext}(E_1,E_2)\, ,
\]
where $E_i = E(S_i,L_i)$.
\end{thm}

\begin{expl}
 Let $S = S_0$ be the zero section of
 $\pi : M = T^*B/\Lambda \to B$ with trivial bundle.
 This Lagrangian submanifold corresponds to the structure sheaf
 $E(S_0) = \mathcal{O}_W$ on $W$.
 A principal polarization (an ample line bundle of degree 1) 
 $E \to W$ corresponds to an affine
 Lagrangian torus $S_1$ of ``slope 1''.
 Similarly, $E^k \to W$ corresponds to an affine Lagrangian torus $S_k$
 of ``slope k''.
 Hence we have
 \[
   \dim H^0(W,E^k) = k^n= \# S_0 \cap S_k = \dim HF(S_0,S_k)\, .
 \]
 These Lagrangian sections $S_k$ are just the section $\lambda_k$
 appeared in the proof of Theorem \ref{thm-Anderson}.
 In other words, $HF(S_0,S_k)$ can be identified with the space of
 wave functions defined by the real polarization
 $\pi : M \to B$.
 In this case, the isomorphism of the spaces of wave functions
 can be regarded as a part of mirror symmetry.

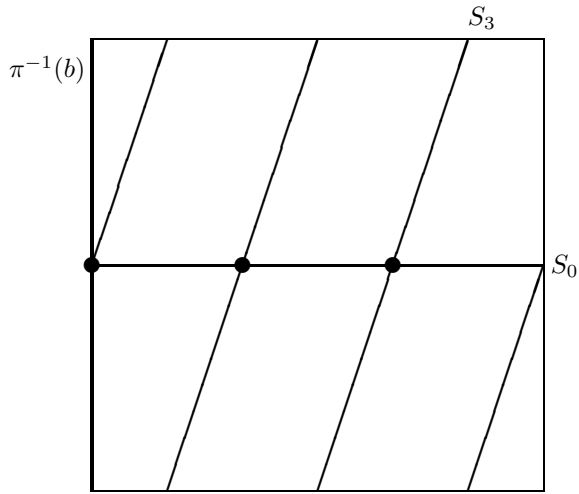
\begin{center}
\setlength{\unitlength}{1mm}
\begin{picture}(70,90)(20,0)
  \put(30,0){\framebox(60,60)}
  \thicklines
  \multiput(40,0)(20,0){2}{\line(1,3){20}}
  \put(30,30){\line(1,0){60}}
  \put(30,30){\line(1,3){10}}
  \put(80,0){\line(1,3){10}}
  \put(30,0){\line(0,1){60}}
  \multiput(30,30)(20,0){3}{\circle*{2}}
  \put(91,29){$S_0$}
  \put(80,62){$S_3$}
  \put(19,55){$\pi^{-1}(b)$}
\end{picture}
\end{center}

\end{expl}

The next theorem states the correspondence of the product structure
(compositions) of morphisms.
\begin{thm}[Fukaya \cite{F}]
The following diagram commutes.
\[
\begin{CD}
  HF \bigl( (S_1,L_1), (S_2,L_2) \bigr)
      \otimes
  HF \bigl( (S_2,L_2), (S_3,L_3) \bigr)
   @>>> 
  HF \bigl( (S_1,L_1), (S_3,L_3) \bigr)\\ 
  @VV{\mathfrak{m}_2}V  @VV{\mathfrak{m}_2}V \\
    \mathrm{Ext}(E_1,E_2)
    \otimes
    \mathrm{Ext}(E_2,E_3)
  @>>> 
  \mathrm{Ext}(E_1,E_3)
\end{CD}
\]
\end{thm}

The matrix elements are given by {\it theta functions}.
We see this in the following example.
\newcommand{\thta}[3]{\vartheta
    \genfrac[]{0pt}{}{#1}{0} (#2 \tau, #3)}

\begin{expl}[Polishchuk-Zaslow \cite{PZ}]
Let $M = T^2 = \mathbb{R}^2/\mathbb{Z}^2$
be a 2-torus 
with a complexified K\"ahler form
$\omega^{\mathbb{C}} = - \sqrt{-1} \tau dx \wedge dy$.
In this case, its mirror is given by
$W = \mathbb{C}/\tau \mathbb{Z} + \mathbb{Z}$.
We consider the affine Lagrangian submanifolds $S_k$ of $M$
corresponding to the lines $y = kx$ in $\mathbb{R}^2$ for $k = 0,1,2$.
Then 
\begin{align*}
  S_0 \cap S_1 &= S_1 \cap S_2 = \{b_0\}\, , \\
  S_0 \cap S_2 &= \{ b_0, b_1\}\, , 
\end{align*}
where $b_0$ and $b_1$ are the points 
corresponding to $(0,0), (\frac 12, 0) \in \mathbb{R}^2$ respectively.

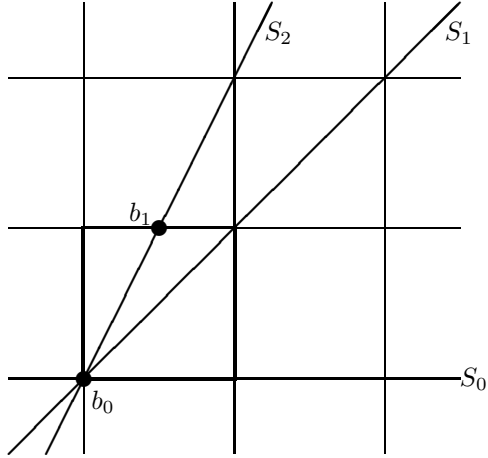
\begin{center}
\setlength{\unitlength}{1mm}
\begin{picture}(70,60)
  \multiput(0,10)(0,20){3}{\line(1,0){60}}
  \multiput(10,0)(20,0){3}{\line(0,1){60}}
  \put(10,10){\circle*{2}}
  \put(20,30){\circle*{2}}
  \thicklines
  \put(10,10){\framebox(20,20)}
  \put(0,10){\line(1,0){60}}
  \put(0,0){\line(1,1){60}}
  \put(5,0){\line(1,2){30}}
  \put(60,9){$S_0$}
  \put(58,55){$S_1$}
  \put(34,55){$S_2$}
  \put(11,6){$b_0$}
  \put(16,31){$b_1$}
\end{picture}
\end{center}
Then
\begin{align*}
   HF(S_0,S_1) &= \mathbb{C} b_0\, ,\\
   HF(S_1,S_2) &= \mathbb{C} b_0\, ,\\
   HF(S_0,S_2) &= \mathbb{C} b_0 \oplus \mathbb{C} b_1
\end{align*}
and the product
\[
  \mathfrak{m}_2 : HF(S_0,S_1) \otimes HF(S_1,S_2) \to HF(S_0,S_2)
\]
is defined by counting triangles bounded by $S_0, S_1, S_2$.
More precisely, 
\[
  \mathfrak{m}_2(b_0 \otimes b_0) =
  C(b_0,b_0,b_0) \cdot b_0 + C(b_0,b_0,b_1) \cdot b_1
\]
with
\[
  C(b_0,b_0,b_0) = \sum \exp (2\pi \sqrt{-1} 
  (\text{area of triangle}))\, ,
\]
where the summation is taken over all triangles with vertices
$b_0,b_0,b_0$, etc.
Consequently we have
\[
  \mathfrak{m}_2(b_0 \otimes b_0) =
  \thta{0}{2}{0} \cdot b_0 + \thta{\frac 12}{2}{0} \cdot b_1\, .
\]

On the other hand, 
$S_k$ corresponds to a holmorphic line bundle on $W$ of degree $k$:
\begin{align*}
  &E_0 = \mathcal{O}_W\, ,\\
  &E_1 = E \,: \text{ a principal polarization}\, ,\\
  &E_2 = E^2.
\end{align*}
Then
\begin{align*}
   \mathrm{Ext}(E_0,E_1) &= H^0(W,E) 
   = \mathbb{C} \cdot \thta{0}{2}{z} \, ,\\
   \mathrm{Ext}(E_1,E_2) &= H^0(W,E)\, ,\\
   \mathrm{Ext}(E_0,E_2) &= H^0(W,E^2) 
    = \mathbb{C} \cdot \thta{0}{2}{2z}  
   \oplus \mathbb{C} \cdot \thta{\frac 12}{2}{2z}
\end{align*}

In this case, the product
\[
  H^0(E) \otimes H^0(E) \to H^0(E^2)
\]
is the natural product
\begin{multline*}
  \thta{0}{2}{z} \thta{0}{2}{z} \\
  = \thta{0}{2}{0} \thta{0}{2}{2z} 
  + \thta{\frac 12}{2}{0} \thta{\frac 12}{2}{2z}
\end{multline*}
(``the addition formula'').
\end{expl}

\begin{rem}
For more general Calabi-Yau manifolds with (special) 
Lagrangian fibrations satisfying suitable conditions,
and its ``topological mirror'' $W \to B$,
we can generalize some parts of the above discussion.

\begin{thm}[Gross \cite{G1}, Tyurin \cite{Ty1}]
Let $M$, $W$ be a mirror pair of K3 surfaces 
and take an ample line bundle $E$ on  $W$.
Let $S_0$ and $S$ be Lagrangian section in $M$ corresponding to
$\mathcal{O}_W$ and $E$ respectively.
Then
\begin{equation}
  (-1)^{n(n-1)/2} S_0 \cdot S = \chi (E)\, .
 \label{mirrorRR}
\end{equation}
\end{thm}

It is not known whether the Floer homology $HF(S_0,S)$ can be defined,
and $\dim HF(S_0,S)$ coincides with the number of the intersection
$S_0 \cap S$.

For higher dimensional case,
(\ref{mirrorRR}) holds asymptotically:
\begin{thm}[Gross \cite{G1}]
Under some assumptions,
\[
  \text{ the leading term of }
   (-1)^{n(n-1)/2} S_0 \cdot S_k 
 = \text{the leading term of } \chi (E^k).
\]
\end{thm}

\end{rem}


\newpage
\section{Projective embeddings and K\"ahler metrics}

By definition, every projective manifold $X$ can be 
embedded into projective spaces.
A way to equip a K\"ahler metric on $X$ is to restrict the 
Fubini-Study metric on the projective space.
On the other hand, some projective manifolds have standard K\"ahler
metrics such as Ricci-flat metrics and K\"ahler-Einstein metrics.
The restriction of the Fubini-Study metric is not a desired one
in general.
For example, the restriction of the Fubini-Study metric to cubic
curves in $\mathbb{CP}^2$ is not flat.
K\"ahler metrics obtained from the Fubini-Study metric is
restricted.

A way to deal with more general K\"ahler metrics 
from the view point of projective embeddings is 
to consider all tensor power of the polarization $L \to X$.
We can think of this as an analogy of the approximation of smooth
functions by polynomials.
In this section, we recall two theories of approximations of
K\"ahler metrics.

\subsection{Bergman kernels and projective embeddings}

Let $X$ be a smooth projective variety
and $L \to X$ an ample line bundle.
In this subsection, we fix a K\"ahler metric $\omega$
in the class $c_1(L)$
and a Hermitian metric on $L$ whose first Chern form coincides with
$\omega$.

For each $k$, we take a basis 
$s_0, \dots , s_{N_k}$ of the space $H^0(X, L^k)$ of 
holomorphic sections, and consider the embeddings defined by
\[
  \iota_k : X \hookrightarrow \cp{N_k},
  \quad z \mapsto (s_0(z): \dots : s_{N_k}(z))\, . 
\]
We denote the pull-back of the Fubini-Study metric by $\omega_k$:
\[
  \omega_k := \frac 1k \iota_k^* \fs\, ,
\]
here we normalize the metric so that
$\omega_k$ represents $c_1(L)$.
Tian and Zelditch proved the following theorem:
\begin{thm}[Tian \cite{T}, Zelditch \cite{Zel}]
  \label{thm-TZ}
  If $s_0, \dots , s_{N_k}$ are orthonormal for each $k \gg 1$,
  then $\omega_k$ converges to the original metric $\omega$
  as $k \to \infty$ in the $C^{\infty}$-topology.
More precisely, 
\begin{equation}
  \| \omega - \omega_k \|_{C^r} = O(k^{-1})
  \label{th-Z}
\end{equation}
for each $r$.
\end{thm}

Namely, any K\"ahler metric representing $c_1(L)$
can be approximated by pull-back of Fubini-Study metrics on
$\cp{N_k} = \mathbb{P} H^0(X,L^k)$ for large $k$.

\begin{rem}
  The case of abelian varieties is discussed by Ji \cite{Ji}
  and Kempf \cite{K}.
\end{rem}

This theorem follows from the asymptotic behavior of the 
Bergman kernels $\Pi_k(z,w)$.
Recall that $\Pi_k(z,w)$ is given by
\[
  \Pi_k(z,w) = \sum_{i=0}^{N_k} s_i(z) s_i(w)^*
\]
for orthonormal basis $s_0, \dots , s_{N_k}$ of $H^0(X,L^k)$.
Its diagonal part
\[
  f_k(z) =\Pi_k(z,z) = \sum_{i=0}^{N_k} | s_i(z) |_h^2
\]
gives a distortion function, i.e.
\[
  \omega - \omega_k = \frac {\sqrt{-1}}k \dd \log f_k \, .
\]
Theorem \ref{thm-TZ} follows from the following theorem.
\begin{thm}[Zelditch \cite{Zel}]\label{expansion}
$f_k$ has the following asymptotic expansion:
\[
  f_k(z) = k^n + a_1 k^{n-1} + a_2 k^{n-2} + \dots
\]
with $a_i(z) \in C^{\infty}(X)$ and
\[
  \left\| f_k - \sum_{i=0}^q a_i k^{n-i} \right\|_{C^r}
  \le C_{r,q} k^{n-q-1}
\]
for some constant $C_{r,q} > 0$,
here we put $a_0 = 1$.
\end{thm}

\subsection{Stability and K\"ahler metrics}

The notion of stability is introduced by Mumford 
to construct moduli spaces.
Moduli spaces are constructed as quotient spaces.
In general, quotient spaces are not Hausdorff.
To obtain a Hausdorff space, we consider moduli space of 
(semi-)stable objects.

On the other hand, stability is closely related to existences of 
differential geometric structures such as metrics and connections.
The following theorem is a typical example,
known as ``Hitchin-Kobayashi correspondence''.
\begin{thm}[Narasimhan-Seshadri \cite{NS}, Donaldson \cite{D1, D2},
  Uhlenbeck-Yau \cite{UY}]
  Let $(X,L)$ be a polarized manifold and $E \to X$ a holomorphic
  vector bundle.
  Then $E$ is stable in the sense of Mumford-Takemoto
  if and only if $E$ admits an irreducible 
  Hermitian-Yang-Mills connection.
\end{thm}

It is expected that such a theorem also holds for the case of manifolds.

\begin{conj}[Hitchin-Kobayashi correspondence for manifolds]
 A smooth polarized projective variety $(X,L)$ admits a K\"ahler metric
 of constant scalar curvature in the class $c_1(L)$ if and only if
 it is stable in a suitable sense.
\end{conj}

There are several notion of stability (these are closely related each
other but not equivalent) and it is not clear which stability is the right
one at the moment.

Recently, Donaldson proved a part of the above conjecture.
Here we review the result
.
\begin{df}
  For a basis $s_0,\ldots,s_{N_k}$ of $H^0(X,L^k)$, put
  \[
    M_{ij} 
     = \sqrt{-1} \int_{\iota_k (X)}
       \frac {Z^i \bar{Z^j}}{\sum_l |Z^l|^2} {\omega_{\mathrm{FS}}}^n
     = \sqrt{-1} \int_X 
       \frac {s_i \overline{s_j}}{\sum_l |s_l|^2} (k \omega_k)^n
  \]
  and $M = (M_{ij}) \in \mathfrak{su}(N_k +1)$.
  The basis
  $s_0,\ldots,s_{N_k}$ is said to be {\it balanced} 
  if $M$ is a constant multiple of the identity matrix.
\end{df}

Since $\omega_k$ changes as $s_0,\ldots,s_{N_k}$,
it is not clear whether balanced basis exist or not.
The existence of balanced basis is closely related to the stability
of the polarized manifold $(X,L)$.
\begin{thm}[Luo \cite{Luo}, Zhang \cite{Zh}]
Assume that $\mathrm{Aut}(X,L)/\mathbb{C}^*$ is discrete.
If $H^0(L^k)$ has a balanced basis,
then the Hilbert/Chow point of $X$ is stable.
\end{thm}

\begin{thm}[Donaldson \cite{D}] \label{thm-D}
  Assume that $X$ admits a K\"ahler metric $\omega_{\infty}$ 
  in the class $c_1(L)$ of constant scalar curvature.
  Further we assume that the automorphism group
  $\mathrm{Aut} (X,L)$ of the polarized variety $(X,L)$
  is discrete.
  Then, for each $k \gg 1$, we can take a {\it balanced} embeddings.
  In particular, $(X,L)$ is stable.
  Furthermore, the pull-backs of the Fubini-Study metrics 
  $\omega_k = \frac 1k \fs$ converges to
  the constant scalar curvature metric $\omega_{\infty}$.
\end{thm}


\newpage
\section{Abelian varieties and theta functions}

In section 3 and 4, we saw some relations between theta functions and 
Lagrangian fibrations on abelian varieties
from the point of view of geometric quantization and mirror symmetry.
It is natural to ask what happen if we embed the abelian varieties
into projective spaces by these theta functions.
In this section, we prove that the Lagrangian fibration of the abelian
variety can be approximated by the restrictions of the moment maps
of projective spaces in a certain sense.


\subsection{Statement of the main theorem}

We consider the same situation as in section 3.5.
Let $X= \mathbb C^n /( \Omega \mathbb Z^n + \mathbb Z^n)$ be 
an abelian variety of complex dimension $n$,
where $\Omega$ is an $n \times n$ complex symmetric matrix whose 
imaginary part $\mathrm{Im}\,\Omega$ is positive definite,
and take an ample symmetric line bundle $L$ on $X$ of degree 1
defined by
\begin{equation}
  L = ( \mathbb C^n \times \mathbb C )/(\Omega \mathbb Z^n + \mathbb Z^n)\, ,
\end{equation}
\[
  (z,\zeta) \sim (z+\lambda, 
           e^{\pi {}^t \lambda (\mathrm{Im}\,\Omega)^{-1}z
               + \frac{\pi}2 {}^t \lambda (\mathrm{Im}\,\Omega)^{-1} \lambda}
            \zeta),
   \quad \lambda \in \Omega \mathbb Z^n + \mathbb Z^n
\]
We fix a Hermitian metric
\[
  h_0 = \exp(-\pi {}^t z (\mathrm{Im}\,\Omega)^{-1} \bar z)
\]
on $L$ and a flat metric
\[
  \omega_0 = 
  c_1(L,h_0) = - \frac{\sqrt{-1}}{2 \pi} \partial \bar\partial \log h_0
             =  \frac{\sqrt{-1}}2 \sum_{i,j} h_{ij} dz^i \wedge d \bar z^j\, .
\]
on $X$,
where we denote $(h_{ij}) = (\mathrm{Im}\,\Omega)^{-1}$.

We consider the following Lagrangian torus fibration as before:
\[
  \pi : X = X^+ \times X^- \longrightarrow X^-= T^n, \quad
  z = \Omega x + y \longmapsto y\, .
\]

For each $k$, we denote $\{s_1,\ldots,s_{k^n}\}$ the basis of 
$H^0(X,L^k)$ defined in (\ref{theta})
and consider the embedding
\[
  \iota_k : X \longrightarrow \mathbb {CP}^{N_k}, \qquad
  x \mapsto (s_0(x):\cdots :s_{N_k}(x))\, .
\]
Write $\omega_k = \frac 1k \iota_k^* \omega_{\mathrm FS}$.

We compare the Lagrangian fibration $\pi : X \to X^-$ 
and the restriction of a natural
Lagrangian fibration (the moment map of the natural torus action) 
on projective spaces.

Let $T^{k^n-1} \subset SU(k^n)$ be the maximal torus which consists of 
diagonal matrices and consider its action on $\mathbb {CP}^{k^n-1}$ given by
\[
  (Z_1:\cdots :Z^{k^n}) \mapsto (\lambda_1 Z^1:\cdots :\lambda_{k^n} Z^{k^n})\, ,
  \quad  \begin{pmatrix}
            \lambda_1 & & 0\\
              & \ddots & \\
            0 & & \lambda_{k^n}
          \end{pmatrix} \in T^{k^n-1}.
\]
Then this action is Hamiltonian and its moment map is given by
\[
  \mu_{T^{k^n-1}}(Z^1:\cdots :Z^{k^n})= \frac 1{\sum |Z^i|^2}
        \left(|Z^1|^2, \dots , |Z^{k^n}|^2 \right)\, ,
\]
here we identify the dual of the Lie algebra of $T^{k^n-1}$ with
\[
  \left\{(\xi_1, \dots , \xi_{k^n}) \in {\mathbb R}^{k^n} \, 
       \Bigm| \, \sum \xi_i = 1 \right\}.
\]
Denote the image of $\iota_k (X) \subset \mathbb {CP}^N$ under $\mu_{T^{k^n-1}}$
by $B_k$ and consider the following map
\[
  \pi_k := \mu_{T^{k^n-1}} \circ \iota_k : X \longrightarrow B_k \, .
\]
Note that $\pi_k$ has many degenerate fibers.
In fact the most degenerate one is 0-dimensional.
We also remark that each fiber is isotropic with respect to
$\omega_k$ (at smooth points).

Next we introduce a distance on $B_k$.
For that purpose, we define a distance on the polytope
$\Delta_k = \mu_{T^{k^n-1}}(\mathbb {CP}^{k^n-1})$.
Note that $\mu_{T^{k^n-1}}$ has no degenerate
fiber on the interior $\overset{\circ}{\Delta}_k$ of $\Delta_k$.
We take a metric on $\overset{\circ}{\Delta}_k$ so that
the moment map 
$\mu_{T^{k^n-1}} : (\mathbb {CP}^{k^n-1}, \omega_{\mathrm{FS}}) \to \Delta_k$
is a Riemannian submersion on the interior.
This is equivalent to the following definition.
Consider the restriction $\mu_{T^{k^n-1}} : \mathbb {RP}^{k^n-1} \to \Delta_k$
of the moment map to $\mathbb {RP}^{k^n-1} \subset \mathbb {CP}^{k^n-1}$.
This is a $2^{k^n-1}$-sheeted covering which branches on the boundary of 
$\Delta_k$.
By identifying $\overset{\circ}{\Delta}_k$ with a sheet of $\mathbb {RP}^{k^n-1}$, 
we have the restriction of the Fubini-Study metric on $\overset{\circ}{\Delta}_k$.
We define a distance on $B_k$ induced by the normalized metric on 
$\overset{\circ}{\Delta}_k$ with normalized factor $\frac 1k$.

On the other hand, We equip the metric on $X^-$ so that 
$\pi : (X, \omega) \to X^-$ is a Riemannian submersion.

Under the above preparation, we can state the main theorem
\begin{thm} \label{main}
  \begin{enumerate}
     \item $\omega_k$ converges to $\omega_0$ in $C^{\infty}$.
           In particular, the sequence $(X,\omega_k)$ 
           of compact Riemanian manifolds converges to $(X,\omega_0)$ 
           in Gromov-Hausdorff topology.
    \item $B_k$ converge to $X^-$ in Gromov-Hausdorff topology.
    \item The sequence of maps $\pi_k : X \to B_k$ between metric spaces 
          converge to $\pi : X \to X^-$.
  \end{enumerate}
\end{thm}

Before the proof of this theorem, 
we recall the definition of Gromov-Hausdorff distances.
If $Z$ is a metric space and $X,Y \subset Z$,
then the {\it Hausdorff distance} of $X$ and $Y$ is defined by
\[
  d^Z_{\mathrm{H}}(X,Y) = \inf \{ \varepsilon > 0 \,\, 
              |\,\, X \subset B(Y,\varepsilon) \,\, 
              \text{and}\,\, Y \subset B(X,\varepsilon)\},
\]
where we denote $B(X,\varepsilon)$ the $\varepsilon$-neighborhood of 
$X$ in $Z$.

For metric spaces $X$ and $Y$, the {\it Gromov-Hausdorff distance} is 
defined by
\[
  \dgh(X,Y) = \inf \{ d^Z_{\mathrm{H}}(X,Y) \,\, 
                 | \,\, X, Y \hookrightarrow Z 
                   \,\,\text{are isometric embeddings.}\}\, .
\] 
This is equivalent to the following definition:
\[
  \dgh(X,Y) = \inf \bigl\{ d^{X \coprod Y}_{\mathrm{H}}(X,Y) \bigr\}\, ,
\]
where the infimum is taken over all metrics on $X \coprod Y$ compatible with 
these on $X$ and $Y$.

Next we define convergence of maps.
Let $X_i,Y_i,X,Y$ be metric spaces and consider the maps
$f_i : X_i \to Y_i$, $f : X \to Y$.
Suppose that $X_i$ and $Y_i$ converge to $X$ and $Y$ respectively
in Gromov-Hausdorff topology.
Then, from the definition of Gromov-Hausdorff distance, 
there exist metrics on $X \coprod (\coprod_i X_i)$  
such that $X_i$ converge to $X$ in Hausdorff topology
in $X \coprod (\coprod_i X_i)$
(and the same is true for $Y_i$ and $Y$). 
In this case, we say that $\{f_i\}$ converges to $f$ if 
for every sequence $x_i \in X_i$ converging to $x \in X$,
$f_i(x_i)$ converges to $f(x)$ in $Y \coprod (\coprod_i Y_i)$. 

\subsection{Proof of Theorem \ref{main}}

Let $\{s_1,\ldots,s_{k^n}\}$ be the basis of $H^0(X,L^k)$ as above.

\begin{prop} \label{onb}
  The basis $\{s_1,\ldots,s_{k^n}\}$ of $H^0(X,L^k)$ are both balanced and 
  orthonormal with respect to $h_0$ and $\omega_0$.
\end{prop}

\begin{proof}
  Recall that the $G_k$-action on $H^0(X,L^k)$ preserves the
  $L^2$-inner product defined by $\omega_0$ and $h_0$.
  From the formula (\ref{rep}),
  \[
    \| s_{b_i} \|^2_{L^2} = \| \rho_k (-b_i) s_{b_1} \|^2_{L^2}
    = \| s_{b_1} \|^2_{L^2}\, ,
  \]
  i.e. all $s_i$'s have the same $L^2$-norm, here we identify $b_i$
  with $ (1,0,b_i) \in G_k$. 
  Furthermore, we have
  \begin{align*}
    e^{2 \pi \sqrt{-1} {}^ta b_i } (s_{b_i}, s_{b_j})_{L^2}
    &= ( \rho_k(a) s_{b_i}, s_{b_j})_{L^2}\\
    &= ( s_{b_i}, \rho_k(-a) s_{b_j})_{L^2}\\
    &= e^{2 \pi \sqrt{-1} {}^ta b_j } (s_{b_i}, s_{b_j})_{L^2}
  \end{align*}
  for all $a \in X_k^+$.
  This implies that $(s_{b_i}, s_{b_j})_{L^2} = 0$ if $i \ne j$.
  In other words, $s_1, \dots , s_{k^n}$ are orthonormal basis up to
  a constant.

  To show that each $s_i$ has unit norm, we consider the following function
  \[
    f_k(z) := \sum_{i=1}^{k^n} |s_i(z)|^2_{h_0}\, .
  \]
  From the definition of $h_0$ and $s_i$,
  \begin{align*}
    |s_i(z)|^2_{h_0}
    &= \exp (-k \pi \xt z (\im \Omega)^{-1} \bar{z} ) \\
    &\quad \cdot
       C_{\Omega}^2 k^{-\frac n2}
       \exp \left( \frac {\pi k}2 \xt z (\im \Omega)^{-1} z \right)
       \exp \left( \frac {\pi k}2 \xt \bar{z}
       (\im \Omega)^{-1} \bar{z} \right)\\
    & \quad \cdot \left| \tht{-b_i}z \right|^2\\
    &= C_{\Omega}^2 k^{-\frac n2}
       \exp \left( \frac {\pi k}2 \xt (z-\bar{z}) (\im \Omega)^{-1}
       (z-\bar{z}) \right)
       \left| \tht{-b_i}z \right|^2\\
    &= C_{\Omega}^2 k^{-\frac n2}
       \exp \left( -2 \pi k \xt (\im z) (\im \Omega)^{-1}
       (\im z) \right)
       \left| \tht{-b_i}z \right|^2.
  \end{align*}
  By using $z = \Omega x + y$,
  we have
  \begin{equation}
   \begin{split}
    |s_i(z)|^2_{h_0} &= C_{\Omega}^2 k^{-\frac n2}
       \exp \left( -2 \pi k \xt x (\im \Omega) x \right)
       \left| \tht{-b_i}z \right|^2\\
    &= C_{\Omega}^2 k^{-\frac n2}
     e \left( \frac k2 \xt x (\Omega - \bar{\Omega}) x \right)
     \left| \tht{-b_i}z \right|^2.
   \end{split}\label{C^0-norm}
  \end{equation}
  From the definition of theta functions,
  \begin{align*}
    &\left| \tht{-b_i}z \right|^2 \\
    &\quad
    = \sum_{l,m \in \mathbb{Z}^n}
       e \left( \frac 1{2k} \xt l \Omega l +\xt l (z - b_i) \right)
       \cdot e \left( -\frac 1{2k} \xt m \bar{\Omega} m
       - \xt m (\bar{z} - b_i) \right)\\
    &\quad
    = \sum_{l,m \in \mathbb{Z}^n}
       e \left( \frac 1{2k} \xt l \Omega l
       -\frac 1{2k} \xt m \bar{\Omega} m
       +\xt l \Omega x
       - \xt m \bar{\Omega} x
       + \xt (l-m) (y - b_i) \right).
  \end{align*}
  Therefore we have
  \begin{align*}
    &|s_i(z)|^2_{h_0} \\
    &\quad = C_{\Omega}^2 k^{-\frac n2}
     \sum_{l,m \in \mathbb{Z}^n}
     e \biggl( \frac k2 (\xt x \Omega x- \xt x \bar{\Omega} x)
     + \frac 1{2k} \xt l \Omega l
     -\frac 1{2k} \xt m \bar{\Omega} m \biggr.\\
    &\phantom{\quad = C_{\Omega}^2 k^{-\frac n2}
     \sum_{l,m \in \mathbb{Z}^n}
     e \biggl( \frac k2 \xt x \Omega x- \xt x \bar{\Omega} x}
    \biggl. +\xt l \Omega x - \xt m \bar{\Omega} x
     + \xt (l-m) (y - b_i) \biggr)\\
    &\quad = C_{\Omega}^2 k^{-\frac n2}
     \sum_{l,m \in \mathbb{Z}^n}
     e \biggl( \frac k2 \lt \left( x+ \frac lk \right) \Omega
     \left( x+ \frac lk \right) - \frac k2
     \mt \left( x+ \frac mk \right) \bar{\Omega}
     \left( x+ \frac mk \right)\\
   &\phantom{\quad = C_{\Omega}^2 k^{-\frac n2}
     \sum_{l,m \in \mathbb{Z}^n}
       e \biggl( \frac k2 \xt x \Omega x- \xt x \bar{\Omega} x
       + \frac 1{2k} \xt l \Omega l }
    +  \xt (l-m) y \biggr) \cdot e \bigl( \xt (m-l) b_i \bigr)\, .
  \end{align*}
  From the fact that $\{b_i\} = \frac 1k \mathbb{Z}^n/\mathbb{Z}^n$,
  we have the following identity:  
  \[
    \sum_{i=1}^{k^n} e(\xt l b_i) =
                     \begin{cases}
                        k^n, & l \in k \mathbb{Z}^n,\\
                        0,   & \text{otherwise}.
                     \end{cases}
  \]
  By using this, we obtain
  \begin{multline}
    f_k(z) = C^2_{\Omega} k^{\frac n2} \sum_{%
            \genfrac{}{}{0pt}{}{l,m \in \mathbb Z^n,}{l-m \in k \mathbb Z^n}}
      e \Biggl( \frac k2 \biggl(
      \lt \left( x+ \frac lk \right) \Omega \left( x+ \frac lk \right) \\
      - \mt \left( x+ \frac mk \right) \bar{\Omega} \left( x+ \frac mk \right)
      \biggr) + \xt (l-m)y \Biggr).
  \label{f_k}
  \end{multline}
  Put
  \[
    \| s_i \|^2_{L^2(h_0)} 
    = \int_X |s_i|^2_{h_0} \frac{{\omega_0}^n}{n!}  =  C\, .
  \]
  Then
  \[
    \int_X f_k(z) \frac{{\omega_0}^n}{n!} = C \cdot h^0(X,L^k) = C k^n
  \]
  by definition.
  Since $\omega_0 = -\sum dx^i \wedge dy^i$,
  we denote the volume form ${\omega_0}^n/n!$ by $dx \, dy$.
  From (\ref{f_k}),
 {\allowdisplaybreaks
  \begin{align*}
    \int_X & f_k(z) \frac{{\omega_0}^n}{n!}\\
        &= C_{\Omega}^2 k^{\frac n2} \int_X \sum_{%
              \genfrac{}{}{0pt}{}{l,m \in \mathbb Z^n,}{l-m \in k \mathbb Z^n}}
            e \Biggl( \frac k2 \biggl(
            \lt \left( x+ \frac lk \right) \Omega \left( x+ \frac lk \right) \\*
        & \phantom{C_{\Omega}^2 k^{\frac n2} \int_X 
                       \sum_{%
              \genfrac{}{}{0pt}{}{l,m \in \mathbb Z^n,}{l-m \in k \mathbb Z^n}} eM}
             - \mt \left( x+ \frac mk \right) \bar{\Omega} \left( x+ \frac mk \right)
             \biggr) + \xt (l-m)y \Biggr)
              dx\,dy\\ 
        &= C_{\Omega}^2 k^{\frac n2} \sum_{l \in  \mathbb Z^n} \int_{X^+}
           e \Biggl( \frac k2
             \biggl(\lt \left(x+ \frac lk \right) \Omega \left(x+ \frac lk \right)\\*
        & \phantom{C_{\Omega}^2 k^{\frac n2} \sum_{l \in  \mathbb Z^n} \int_{X^+}
                    e M \frac k2 M \left( x+ \frac lk \right) \Omega }  
             -\lt \left(x+ \frac lk \right) \bar{\Omega} \left(x+ \frac lk \right)
             \biggr) \Biggr) dx\\ 
        &= C_{\Omega}^2 k^{\frac n2} \sum_{l \in  \mathbb Z^n} \int_{X^+}
           e \Biggl( \frac k2
             \left(\lt \left(x+ \frac lk \right) (2 \sqrt{-1} \mathrm{Im}\Omega) 
             \left(x+ \frac lk \right)\right)\Biggr) dx\\
        &= C_{\Omega}^2 k^{\frac n2} \sum_{l \in  \mathbb Z^n} \int_{X^+}
           \exp \left( -2 \pi k 
             \lt \left( x+ \frac lk \right) ( \mathrm{Im}\,\Omega)
             \left(x+ \frac lk \right) \right) dx\\
        &= C_{\Omega}^2 k^{\frac n2} \sum_{%
            \genfrac{}{}{0pt}{}{l \in  \mathbb Z^n,} 
            {m \in \frac 1k \mathbb Z^n/ \mathbb Z^n}} \int_{X^+}
             \exp  \left( -2 \pi k \xt (x+l+m) (\mathrm{Im}\Omega) (x+l+m) \right) dx\\
        &= C_{\Omega}^2 k^{\frac n2} 
                  \sum_{m \in \frac 1k \mathbb Z^n/ \mathbb Z^n} \int_{\mathbb R^n}
             \exp \left( -2 \pi k \xt (x+m) (\mathrm{Im}\Omega) (x+m) \right) dx\\
        &= C_{\Omega}^2 k^{\frac {3n}2} \int_{\mathbb R^n}
             \exp \left( -2 \pi k \xt x (\mathrm{Im}\Omega) x \right) dx\\
        &= C_{\Omega}^2 k^n \sqrt{\det (\mathrm{Im} \, \Omega)^{-1}} 
           \cdot 2^{- \frac n2}
         = k^n\, .
  \end{align*}
  }
  Therefore, $C = \|s_i\|^2_{L^2(h_0)} = 1$.

  Next we prove that these basis satisfies the balanced condition.
  Let $h_k$ be the pull back of the standard Hermitian metric 
  (the Fubini-Study metric) on the hyperplane
  bundle $\mathcal O_{\mathbb{CP}^{k^n-1}}(1)$ {\it i.e.} 
  $h_k = \left( \sum_i |s_i|^2 \right)^{-1}$.
  Then $\omega_k = - \frac{\sqrt{-1}}{2 \pi k} \partial \bar\partial \log h_k$.
  From \eqref{rep} we have 
  \[
    h_k = \left( \frac 1{|\mu_k|} \frac 1{|X^+_k|} 
    \sum_{g \in G_k} | \rho_k(g) s_1 |^2 \right)^{-1},
  \]
  where $\rho_k(g) : H^0(X,L^k) \to H^0(X,L^k)$ is the Heisenberg representation
  of $G_k$.
  This means that $h_k$ and $\omega_k$ are invariant under the
  $G_k$-action.
  Hence this action on $H^0(X,L^k)$ is also unitary
  with respect to $h_k$ and $\omega_k$.
  We can prove that  $s_1, \dots , s_{k^n}$ are balanced by the same argument
  as above.
\end{proof}


The next lemma is a key of the proof of Theorem \ref{main}.
\begin{lem}\label{lem}
  There exists a constant $C>0$ independent of $k$ such that for each
  $z = \Omega x + y \in X$,
  \[
    C^{-1} k^{\frac n2} e^{-kC \mathrm{dist}(y,b_i)^2} \le
      |s_i(z)|_{h_0}^2
      \le C k^{\frac n2} e^{-kC \mathrm{dist}(y,b_i)^2}\, ,
  \]
  for some distance $\mathrm{dist}$ on $X^-$.
\end{lem}

\subsubsection*{Proof.}
For $z = \Omega x + y \in X$, we write 
$w = \Omega^{-1}(z-b_i) = x + \Omega^{-1}(y-b_i)$.
Then
\begin{align*}
 \tht{-b_i}{z}
  &= \sum_{l \in \mathbb{Z}^n}
     e\left( \frac 1{2k} \xt l \Omega l + \xt l(z-b_i) \right)\\
  &= \sum_{l \in \mathbb{Z}^n}
     e\left( \frac 1{2k} \xt l \Omega l +
     \xt l(\Omega x + y-b_i) \right)\\
  &= \sum_{l \in \mathbb{Z}^n}
     e \Biggl( \frac 12 \lt
     \left( \frac l{\sqrt{k}} + \sqrt{k} w \right) \Omega
     \left( \frac l{\sqrt{k}} + \sqrt{k} w \right) \Biggr)\\
  &\quad  \quad \cdot
    e \left( - \frac k2
    \xt \bigl(x + \Omega^{-1}(y-b_i) \bigr) \Omega
    \bigl(x + \Omega^{-1}(y-b_i) \bigr) \right).
\end{align*}
From (\ref{C^0-norm}), we have
\begin{align*}
  |s_i(z)|^2_{h_0} 
  &= C_{\Omega}^2 k^{-\frac n2}
     e \left( \frac k2 \xt x (\Omega - \bar{\Omega}) x \right)\\
  & \qquad \cdot \left| e \left( - \frac k2
    \xt \bigl(x + \Omega^{-1}(y-b_i) \bigr) \Omega
    \bigl(x + \Omega^{-1}(y-b_i) \bigr) \right) \right|^2\\
  & \qquad \cdot \left|
    \sum_{l \in \mathbb{Z}^n}
     e \Biggl( \frac 12 \lt
     \left( \frac l{\sqrt{k}} + \sqrt{k} w \right) \Omega
     \left( \frac l{\sqrt{k}} + \sqrt{k} w \right) \Biggr)
     \right|^2\\
  & = C_{\Omega}^2 k^{-\frac n2}
    \exp \bigl( 2 \pi k \xt (y-b_i) (\im \Omega^{-1}) (y-b_i) \bigr)\\
  & \qquad \cdot \left|
    \sum_{l \in \mathbb{Z}^n}
     e \Biggl( \frac 12 \lt
     \left( \frac l{\sqrt{k}} + \sqrt{k} w \right) \Omega
     \left( \frac l{\sqrt{k}} + \sqrt{k} w \right) \Biggr)
     \right|^2.
\end{align*}
Since
\begin{multline*}
  \frac 1{\sqrt{k}^n} \sum_{l \in \mathbb{Z}^n}
     e \Biggl( \frac 12 \lt
     \left( \frac l{\sqrt{k}} + \sqrt{k} w \right) \Omega
     \left( \frac l{\sqrt{k}} + \sqrt{k} w \right) \Biggr)\\
   \quad =
    \sum_{l \in \frac 1{\sqrt k}\mathbb{Z}^n}
     \frac 1{\sqrt{k}^n}
     e \Biggl( \frac 12 \xt
     \left( l + \sqrt{k} w \right) \Omega
     \left( l + \sqrt{k} w \right) \Biggr)
\end{multline*}
converges to
$\int_{\mathbb{R}^n} e \left( \frac 12 \xt u \Omega u \right) du$
uniformly on a fundamental domain of $\Omega \mathbb{Z}^n + \mathbb{Z}^n$
as $k \to \infty$,
this part can be bounded uniformly from above and below.
Furthermore $\exp \bigl( 2 \pi k \xt (y-b_i) (\im \Omega^{-1}) (y-b_i) \bigr)$
is of the form
$e^{-kC \mathrm{dist}(y,b_i)^2}$,
since $\im \Omega^{-1}$ is negative definite.
Hence we have the desired estimate.
\hfill $\square$

\subsubsection*{Proof of Theorem \ref{main} (1).}
The convergence of K\"ahler forms $\omega_k$ follows from 
Proposition \ref{onb} and Theorem \ref{thm-TZ} or \ref{thm-D}.

\begin{lem}\label{geod}
  Let $d_0$ and $d_k$ denote the distance on $X$ defined by 
  $\omega_0$ and $\omega_k$ respectively.
  Then there exists a constant $C$ independent of $k$ such that
  \[
    \left( 1- \frac Ck \right) d_0(p,q) \le d_k(p,q) \le
       \left( 1+ \frac Ck \right) d_0(p,q)
  \]
  holds for each $p,q \in X$.
\end{lem}

\begin{proof}
  Set $x^{n+i} = y^i$ for $i = 1, \dots ,n$.
  From the estimate \eqref{th-Z} by Zelditch, 
  the equation of geodesics with respect to $\omega_k$ is
  \[
    \ddot x^{\alpha} + \sum_{\beta, \gamma}
       \Gamma^{\alpha}_{\beta \gamma} \dot x^{\beta} \dot x^{\gamma} = 0 \,\,\,
    \text{with} \,\,\, |\Gamma^{\alpha}_{\beta \gamma}| = O(k^{-1}).
  \]

  Let $\gamma_k (t) = \xi_0 t + \eta (t)$, $0 \le t \le 1$
  be a geodesic with respect to $\omega_k$,
  where $\xi_0$ is a constant vector and $\eta (0) = \dot{\eta}(0) = 0$.
  Note that $\Gamma_0(t) = \xi_0 t$ is a geodesic with respect to $\omega_0$.
  Then $\eta$ satisfies
  \[
    \ddot{\eta} + \Gamma (\xi_0 + \dot{\eta}, \xi_0 + \dot{\eta}) =0\, ,
  \]
  where $\Gamma (\dot{x}, \dot{x}) 
  = (\sum_{\beta, \gamma}\Gamma^{\alpha}_{\beta \gamma} 
  \dot x^{\beta} \dot x^{\gamma})$.
  Hence we have
  \[
    \frac d{dt} |\xi_0 + \dot{\eta}|^2
    = 2 (\ddot{\eta}, \xi_0 + \dot{\eta})
    = -2 (\Gamma (\xi_0 + \dot{\eta}, \xi_0 + \dot{\eta}),\xi_0 + \dot{\eta})\, .
  \]
  From the estimate of $\Gamma$,
  \[
    \left| \frac d{dt} |\xi_0 + \dot{\eta}|^2 \right|
    \le \frac Ck |\xi_0 + \dot{\eta}|^3
  \]
  for some constant $C >0$.
  Hereafter we denote constants independent of $k$ be the same $C$.
  Therefore
  \[
    \left| \frac d{dt} \frac 1{|\xi_0 + \dot{\eta}|} \right|
    = \left| \frac 1{|\xi_0 + \dot{\eta}|^2}
     \frac d{dt} |\xi_0 + \dot{\eta}| \right|
     \le  \frac Ck\, .
  \]
  This implies that
  \[
    \left| \frac 1{|\xi_0 + \dot{\eta}|} - \frac 1{|\xi_0|} \right|
    \le \frac {Ct}k \, .
  \]
  From this, we have a uniform bound $|\xi_0 + \dot{\eta}| \le C'$
  for large $k$.
  Hence we have
  \[
    | \ddot{\eta}| = | \Gamma ( \xi_0 + \dot{\eta},\xi_0 + \dot{\eta}) |
    \le \frac Ck \, .
  \]
  Since $\dot{\eta}(0) =0$, we have $|\dot{\eta} | = O \left( \frac 1k \right)$.
  In particular
  \[
    |\dot{\gamma}_k|_{\omega_k} = |\xi_0| + O \left( \frac 1k \right).
  \]
  Therefore the length of $\gamma_k$ is given by
  \[
    \int_0^1 |\dot{\gamma_k}|_{\omega_k} dt
    = |\xi_0| + O \left( \frac 1k \right)
    \int_0^1 |\dot{\gamma_0}|_{\omega_0} dt + O \left( \frac 1k \right).
  \]
  Lemma follows from this estimate.
\end{proof}

Let $\{x_{ij} = \Omega a_i + b_j \}_{i,j = 1,\dots,k^n} = X_k$ 
with the distance induced by $d_0$
and denote the same set with distance $d_k$ by $\{y_{ij}\}$.
Then 
\begin{equation}
 \begin{split}
  \dgh \bigl( (X,\omega_0),\{x_{ij} \} \bigr) &= O(k^{-1})\, ,\\
  \dgh \bigl( (X,\omega_k),\{y_{ij}\} \bigr)  &= O(k^{-1})\, .
 \end{split}
 \label{est1}
\end{equation}
Define a distance on $\{x_{ij}\} \coprod \{y_{ij}\}$ by
\begin{eqnarray*}
  d^{(k)}(x_{ij},y_{ij}) &=& \frac Ck \, ,\\
  d^{(k)}(x_{ij},y_{hl}) &=& 
         \min_{p,q} \left\{ d_0(x_{ij},x_{pq}) + \frac Ck +
                    d_k(y_{pq},y_{hl})\right\}\, ,
\end{eqnarray*}
where $C$ is the constant in Lemma \ref{geod}.
Then 
\[
  \dgh(\{x_{ij} \},\{y_{ij} \}) \le \frac Ck\, .
\]
By combining this with \eqref{est1}, we have
\begin{equation}
  \dgh \bigl( (X,\omega_0),(X,\omega_k) \bigr) = O(k^{-1})\, .
  \label{est2}
\end{equation}
This prove (1) of the theorem.
\hfill $\square$


\subsubsection*{Proof of Theorem \ref{main} (2).\footnote{%
The proof in \cite{N} is not correct.}}
To prove the Gromov-Hausdorff convergence, it suffices to construct
{\it $\ep$-Hausdorff approximations} $\varphi_k : X^- \to B_k$
for large $k$ (see \cite{F2}).

\begin{df}
  Let $(X,d_X)$, $(Y,d_Y)$ be two compact metric spaces.
  A map $\varphi : X \to Y$ is said to be an 
  $\ep$-Hausdorff approximation
  if the following two conditions are satisfied.
  \begin{enumerate}
    \item The $\ep$-neighborhood of $\varphi (X)$ coincides with $Y$.
    \item For each $x,y \in X$,
      \[
        |d_X(x,y) - d_Y(\varphi(x), \varphi(y)) | < \ep \, .
      \]
  \end{enumerate}
\end{df}

We identify $X^-$ with a section $\{0\} \times X^- \subset X$ and put
\[
  \varphi_k := \pi_k|_{X^-} : X^- \longrightarrow B_k \, .
\]
We prove that $\varphi_k$ is a $C/\sqrt{k}$-Hausdorff approximation.

\begin{lem}\label{claim1}
  There exist a constant $C>0$ such that
  the $C/k$-neighborhood of $\varphi_k(X^-)$ is equal to
  $B_k$.
\end{lem}

\begin{proof}
Take any point $\pi_k(z) \in B_k$ and put $y = \pi (z) \in X^-$.
Since $\pi_k : (X,\omega_k) \to B_k$ is invariant under the action of
$X_k^+$, there exists $z' \in \pi^{-1}(y)$ such that
$\pi_k(y) = \pi_k(z')$ and $d_k(z,z') \le C/K$.
Then
\[
  d^{B_k}(\pi_k(z), \varphi_k(y)) = d^{B_k}(\pi_k(z), \pi_k(z'))
  \le d_k(z,z') \le C/K \,.
\]
\end{proof}

For each $\xi \in T_p \cp{N_k}$, we denote its vertical and horizontal
parts by
\[
  \begin{matrix}
    T_p \cp{N_k} & = & T_{\cp{N_k}/\Delta_k , p} & \oplus &
    (T_{\cp{N_k}/\Delta_k , p})^{\perp}\\
    \vin & & \vin & & \vin \\
    \xi & = & \xi^V & + & \xi^H 
  \end{matrix}
\]
where $T_{\cp{N_k}/\Delta_k ,p} = \ker d \mu_k$ is the tangent space
of the fiber of $\mu_k : \cp{N_k} \to \Delta_k$ 
and $(T_{\cp{N_k}/\Delta_k , p})^{\perp}$
is its orthogonal complement with respect to the Fubini-Study metric.
Let $(Z^0: \dots : Z^{N_k})$ be the homogeneous coordinate and write
\[
  \log \frac{Z^i}{Z^0} = u^i + \sqrt{-1} v^i \, .
\]
Then $T_{\cp{N_k}/\Delta_k}$ and $(T_{\cp{N_k}/\Delta_k})^{\perp}$
are spanned by $\frac{\dr}{\dr v^i}$'s and $\frac{\dr}{\dr u^i}$'s
respectively.

Let $\gamma : [0,l] \to B_k$ be a curve and take a lift
$\tilde{\gamma} : [0,l] \to X$ of $\gamma$.
Then the length of $\gamma$ is given by
\[
  \int_0^l \left| \left( \frac d{dt} \tilde{\gamma} 
  \right)^H \right|_{\omega_k} dt
  \, .
\]

\begin{lem}\label{claim2}
\begin{multline*}
  s_j(z) = Ck^{\frac n4} \exp 
  \left( \frac{\pi k}2 \xt z (\im \Omega)^{-1} z \right)
  \exp \left( \frac {\pi k}{\sqrt{-1}} \xt (z-b_j) \Omega^{-1} (z-b_j)
  \right)\\
  \times \left( 1 +  \phi \right)
\end{multline*}
with
\[
  |\phi| = O \left( \frac 1{\sqrt{k}} \right)\,,
  \quad |d \phi | = O(1)\,.
\]
\end{lem}

\begin{proof}
From the definition of the theta function, we have
\begin{align*}
  &\tht{-b}{z} \\
  &\quad = \sum_{l \in \mathbb Z^n} 
       e \left( \frac 12 \xt l \Omega l + \xt l(z+b) \right)\\
  &\quad = \sum_{l \in \mathbb Z^n} 
       e \left( \frac 12 \lt \left( \frac l{\sqrt{k}}
       + \sqrt{k} \Omega^{-1} (z-b) \right) \Omega
       \left( \frac l{\sqrt{k}} + \sqrt{k} \Omega^{-1} (z-b) \right)
       \right)\\
  &\quad \phantom{= \sum_{l \in \mathbb Z^n} 
       e \left( \frac 12 \lt \left( \frac l{\sqrt{k}}
       + \sqrt{k} \Omega^{-1} (z-b) \right) \right)}
    \times e \left( - \frac k2 \xt (z-b)\Omega^{-1} (z-b) \right)
\end{align*}
and
\begin{align*}
 & \sum_{l \in \mathbb Z^n} 
       e \left( \frac 12 \lt \left( \frac l{\sqrt{k}}
       + \sqrt{k} \Omega^{-1} (z-b) \right) \Omega
       \left( \frac l{\sqrt{k}} + \sqrt{k} \Omega^{-1} (z-b) \right)
       \right)\\
 & \quad = k^{\frac n2}
       \sum_{l \in \frac 1{\sqrt{k}} \mathbb{Z}^n + 
       \sqrt{k} \Omega^{-1} (z-b)}
       \exp \left( \pi \sqrt{-1} \xt l \Omega l \right) 
       \frac 1{\sqrt{k^n}}\\
 & \quad = k^{\frac n2} \left( \int_{\mathbb{R}^n} 
       \exp \left( \pi \sqrt{-1} \xt l \Omega l \right)  dl 
       + O \left( \frac 1{\sqrt{k}} \right) \right)\\
 & \quad = Ck^{\frac n2}
  \left(1 + O \left( \frac 1{\sqrt{k}} \right) \right)\, .
\end{align*}
Similarly we have
\[
  \left| d \sum_{l \in \mathbb Z^n} 
       e \left( \frac 12 \lt \left( \frac l{\sqrt{k}}
       + \sqrt{k} \Omega^{-1} (z-b) \right) \Omega
       \left( \frac l{\sqrt{k}} + \sqrt{k} \Omega^{-1} (z-b) \right)
       \right) \right| = O \left( k^{\frac n2} \right)
\]
Lemma \ref{claim2} follows from this.
\end{proof}

From this Lemma we have
\begin{align}
  \frac{Z^j}{Z^0} &= \frac{s_j(z)}{s_0(z)} \notag \\
  &= C_j \exp \Bigl( 2 \pi k \sqrt{-1} \xt (b_j - b_0) \Omega^{-1} z
     \Bigr) + O \left( \frac 1{\sqrt{k}} \right) \notag \\
  &= C_j \exp \Bigl( 2 \pi k \sqrt{-1} \xt (b_j - b_0) 
     (x + \re (\Omega^{-1}) y)
     - 2 \pi k \xt (b_j - b_0) \im (\Omega^{-1}) y
     \Bigr) \notag \\
  & \phantom{= C_j \exp \Bigl( 2 \pi \sqrt{-1} \xt (b_j - b_0) 
     (x + \re (\Omega^{-1}) y)
     - 2 \pi  \xt (b_j - b_0)\Bigr)}
  + O \left( \frac 1{\sqrt{k}} \right)
  \label{non-homog}
\end{align}
for some constant $C_j$.

Let $\gamma : [0,l] \to X^-$ be a curve and take a horizontal lift
$\tilde{\gamma} : [0,l] \to X$, i.e. 
\[
  \frac d{dt} \tilde{\gamma} \in  (T_{X/X^-})^{\perp} \, .
\]
Then the length of $\gamma$ is given by
\[
  \int_0^l \left|  \frac d{dt} \tilde{\gamma} 
  \right|_{\omega} dt
  \, .
\]
Note that 
$(T_{X/X^-})^{\perp}$ is spanned by $J \frac{\dr}{\dr x^i}$'s,
where $J$ is the complex structure on $X$.
(Recall that $\frac{\dr}{\dr x^i}$ tangents to fibers.)
By direct computation, we have
\begin{multline*}
  \left( J \frac{\dr}{\dr x^1} ,\dots , J \frac{\dr}{\dr x^n} \right)
  = \left(\frac{\dr}{\dr x^1} ,\dots ,\frac{\dr}{\dr x^n} \right)
   \Bigl( - \re (\Omega^{-1}) (\im \Omega^{-1})^{-1} \Bigr)\\
   + \left(\frac{\dr}{\dr y^1} ,\dots ,\frac{\dr}{\dr y^n} \right)
   (\im \Omega^{-1})^{-1} 
\end{multline*}
and
\begin{align*}
  \frac{\dr}{\dr x^i} 
    \Bigl( \xt (b_j - b_0) \im (\Omega^{-1}) y \Bigr) = 0,\\
  J \frac{\dr}{\dr x^i}
    \Bigl( \xt (b_j - b_0)(x + \re (\Omega^{-1}) y) \Bigr)= 0\,.
\end{align*}
From this and (\ref{non-homog}), we have
\[
  \left| \frac d{dt} \tilde{\gamma} - 
  \left( \frac d{dt} \tilde{\gamma} \right)^H \right| 
  \le \frac C{\sqrt{k}} 
  \left| \frac d{dt} \tilde{\gamma} \right| \, .
\]
This implies that 
\[
  \Bigl| ( \text{length of $\gamma$}) -
   \text{(length of $\varphi_k(\gamma)$)} \Bigr| 
   \le O \left( \frac 1{\sqrt{k}} \right)\, .
\]
From this, we obtain
\begin{lem}\label{claim3}
For $p,q \in X^-$,
\[
  \big| d^{X^-}(p,q) - d^{B_k}( \varphi_k (p), \varphi_k (q) ) \bigr|
  \le O \left( \frac 1{\sqrt{k}} \right)\, .
\]
\end{lem}

Lemma \ref{claim1} and \ref{claim3} prove the second statement 
of Theorem \ref{main}.
\qed


\subsubsection*{Proof of Theorem \ref{main} (3).}

From the second statement of Theorem \ref{main},
there is a distance on $X^- \coprod B_k$ compatible with those
on $X^-$ and $B_k$ such that
\[
  d^{X^- \coprod B_k}(y, \varphi_k(y)) 
  \le \frac C{\sqrt{k}}
\]
for all $y \in X^-$.

It suffices to show that $\pi_k(p)$ converges to 
$\pi(p)$ for each $p \in X$.
For $p \in X$, 
we denote the corresponding point in $X/X^+_k$ by $\bar p$.
We regard that $\pi(p) \in B \subset X/X^+_k$.
Then the distance between $\pi(p)$ and $\bar p$ 
with respect to $\omega_k$ is 
\[
  d^{X/X^+_k}( \pi(p), \bar p) \le O \left(\frac 1k \right) \,.
\]
Hence
\[
  d^{B_k}(\varphi_k( \pi(p) ) , \pi_k(p))
  \le d^{X/X^+_k}( \pi(p), \bar p) \le O \left(\frac 1k \right) 
\]
and we have
\begin{align*}
  &d^{X^- \coprod B_k}(\pi(p), \pi_k(p))\\ 
  & \quad \le d^{X^- \coprod B_k}(\pi(p), \varphi_k( \pi(p) ))
    + d^{B_k}(\varphi_k( \pi(p) ) , \pi_k(p))\\
  & \quad \le O \left( \frac 1{\sqrt{k}} \right) \,.
\end{align*} 

\qed


\newpage
\section{Lagrangian fibrations and holomorphic sections}

In the previous section, we proved that theta functions can be 
reconstructed from a Lagrangian fibration by using the Bergman kernels.
This construction can be applied to more general settings.
In this section, we study the asymptotic behavior of the
holomorphic sections constructed in this way
by using the asymptotic behavior of the Bergman kernels
(\cite{BSZ}, \cite{Lin}, \cite{B}).


\subsection{Asymptotic behavior of the Bergman kernels}

The asymptotic behavior of the Bergman kernel $\Pi_k(z,w)$ 
is studied in detail.

\begin{thm}[Lindholm \cite{Lin}, Berndtsson \cite{B}]
  \label{Th-LB}
  \[
      \left| \Pi_k(z,w) \right| \le C k^n e^{-c \sqrt{k} \, d(z,w)}
  \]
  for some constants $C, c > 0$,
  where $d(z,w)$ is the distance of $z,w \in X$.
\end{thm}
 
The asymptotic behavior of $\Pi_k(z,w)$  near the diagonal
can be written as follows.
Fix a point $z_0 \in X$ and take a neighborhood $U$ of $z_0$ with
local holomorphic coordinate
$(z^1, \dots ,z^n)$ around $z_0$.
We write 
$\omega = \sqrt{-1} \sum g_{i \bar{j}} dz^i \wedge d\bar{z}^j$
and set $G = ( g_{i \bar{j}} )$.
Then we can take a holomorphic local frame $e_L$ of $L$
satisfying
\[
  \frac {\dr h}{\dr z^i} = \frac {\dr^2 h}{\dr z^i \dr z^j} = 0
  \quad \text{at $z_0$}
\]
for each $i,j = 1, \dots , n$.
Here we put $h(z) = h \left( e_L(z), e_L(z) \right)$.
In other words,
\begin{align*}
  h(z) 
  & = \exp \left( -\sum_{i,j} g_{ij}z^i \bar{z}^j + O( |z|^3) \right)\\
  & = \exp \bigl( - \xt z G \bar{z} + O( |z|^3) \bigr)\, .
\end{align*}
From now on, we fix a trivialization of $L$ over $U$ defined by
\[
  L|_U \longrightarrow U \times \mathbb{C}\, , \quad
  \zeta e_L(z) \longmapsto (z, h(z)^{\frac 12} \zeta)
\]
and express sections of $L^k$ with respect to this trivialization.
Note that this trivialization is not holomorphic but unitary.
Under this setting, $\Pi_k(z,w)$ can be write down as follows:
\begin{thm}[Bleher-Shiffman-Zelditch \cite{BSZ}]
 \begin{align*}
  \Pi_k \left( z_0 + \frac u{\sqrt{k}}, z_0 + \frac v{\sqrt{k}} \right)
  &= \left( \frac k{2 \pi} \right)^n e^{- \frac 12 {}^t uG \bar{u} 
      - \frac 12 {}^t vG \bar{v} +{}^t uG \bar{v}}
     \left( 1 + O \left(\frac 1{\sqrt k} \right) \right)\\
  &= \left( \frac k{2 \pi} \right)^n
     e^{- \frac 12 {}^t (u-v)G \overline{(u-v)}
     + \sqrt{-1} \mathrm{Im}({}^t uGv)}
     \left( 1 + O \left(\frac 1{\sqrt k} \right) \right).
 \end{align*}
\end{thm}

\begin{rem}
This expression is slightly different from that in \cite{BSZ}.
This is caused by the difference of the normalizations of K\"ahler
metrics(in \cite{BSZ}, $\pi c_1(L,h) = \frac 12 \omega$
is used instead of $\omega = 2\pi c_1(L,h) $).
\end{rem}

\subsection{Asymptotic behavior of holomorphic sections}
Let $(L,h) \to (X,\omega)$ be a holomorphic Hermitian line bundle
on a compact K\"ahler manifold such that 
$\frac 1{2\pi} c_1(L,h) = \omega$.

Suppose that we have a Lagrangian fibration
$\pi : (X,\omega) \to B$.
We take $k$-Bohr-Sommerfeld fibers
$\pi^{-1}(b_0) , \dots ,\pi^{-1}(b_N) \subset X$.
Here we consider only the part away from singular fiber.
Thus we assume that $\pi^{-1}(b_0) ,\dots ,  \pi^{-1}(b_N)$ are smooth
(hence $\pi^{-1}(b_i) \cong T^n$).
For each $i$, we take a section
$\sigma_i \in \Gamma (\pi^{-1}(b_i), L^k)$ satisfying 
the following conditions:
\begin{itemize}
  \item $\nabla_{\xi} \sigma_i = 0$ for each 
        $\xi \in T \bigl( \pi^{-1}(b_i) \bigr)$ and
  \item $\int_{\pi^{-1}(b_i)} | \sigma_i |^2_h = 1$.
\end{itemize}
Then 
\[
  \left| \sigma_j (y) \right|^2_h \equiv \frac 1{V_j}\, ,
\]
where $V_i$ is the volume of the fiber $\pi^{-1}(b_j)$.
Put
\[
  s_i(z) = \left( \frac k{2 \pi} \right)^{- \frac n4}
           \int_{\pi^{-1}(b_i)} \Pi_k(z,y) \sigma_i(y) dvol\, ,
\]
where $dvol \in \Omega^n(X)$ is the volume form of fibers
defined by the K\"ahler metric on $X$.

We study the asymptotic behavior of $s_i$'s by using the results
quoted in the previous subsection.
All the proofs of the results in this subsection are given in the next
subsection.

Hereafter we fix an arbitrary constant $\ep >0$.
Let $\rho >0$ be a constant such that local holomorphic coordinates
considered in the previous subsection 
can be defined in the ball of radius $\rho$.
We take a constant $R >0$ such that 
\[
  \min \left\{
  e^{- \frac c{\sqrt{2}} R} , \, e^{- \frac R2} , \,
  \frac {\rho}R \right\} \le \ep\, .
\]
Furthermore we take $r >0$ such that for each $z \in X$ and
$b \in B$ with $d( \pi (z), b) \le r$,
there exists $y_0 \in \pi^{-1}(b)$ with
$ d(z, y_0) = d(z, \pi^{-1}(b))$ and
\[
  d(z,y) \ge \frac 12 d(y_0, y)
\]
for all $y \in \pi^{-1}(b)$.
Let $k$ be large enough so that 
\[
  \frac R{\sqrt{k}} \le \rho \quad \text{and} \quad
  \frac 1{\sqrt{k}} \left( 3R + \frac{3n}{4c} \log k \right) \le r\, .
\]

\begin{prop} \label{decay}
\begin{enumerate}
  \item If $d(z,\pi^{-1}(b_j)) \le \frac 1{\sqrt{k}}
        \left( 3R + \frac{3n}{4c} \log k \right)$,
        then
   \[
     |s_j(z)|_h \le C k^{\frac n4} 
     e^{- \frac c4 \sqrt{k} d(z,\pi^{-1}(b_j))}\, .
   \]
  \item If $d(z,\pi^{-1}(b_j)) \ge \frac 1{\sqrt{k}}
        \left( 3R + \frac{3n}{4c} \log k \right)$,
        then
    \[
       |s_j(z)|_h = O(\ep)\, .
    \]
\end{enumerate}
\end{prop}

Next we give the behavior of $s_j$ near $\pi^{-1}(b_j)$.

Fix $z_0 \in \pi^{-1}(b_0)$ and take an action-angle coordinate
$(x^1, \dots , x^n,y^1, \dots , y^n)$ around $\pi^{-1}(\pi (z_0))$ 
constructed in section 2.
Then we can take a holomorphic coordinate around $z_0$
of the form $z = y + \Omega x + O(|z|^2)$  
with some symmetric matrix $\Omega$ with $\mathrm{Im} \, \Omega > 0$.
In this case, $\omega = \sqrt{-1} \sum g_{i \bar{j}} 
dz^i \wedge d\bar{z}^j = - \sum dx^i \wedge dy^i$ and
$\mathrm{Im} \, \Omega = \frac 12 G^{-1}$ at $z_0$.
Note that $b_j \in \frac 1k \mathbb{Z}^n$ with respect to this 
action angle coordinate.

\begin{thm} \label{behavior}
  For $z \in X$ with $d(z,z_0) \le R/ \sqrt k$ and
  $b=b_j$ with $d(b_0, b_j) \le R/ \sqrt k$,
  \begin{multline*}
     s_j(z) = 
     \left( \frac k{2 \pi} \right)^{\frac n4} \frac 1{\sqrt{V_j}}
     \exp \left( \frac k2 \xt zG (z - \bar{z}) 
     - \sqrt{-1}k{}^tzb - \frac k4 {}^tb(\mathrm{Im}\, \Omega)b \right)\\
     \cdot \left( 1+ O \left( \frac{1}{\sqrt{k}} \right)
     + O(\ep) \right).
  \end{multline*}
\end{thm}

By using the above results, 
we compute $L^2$-inner products of $s_i$'s.

\begin{prop}\label{inner-prod-1}
  If $d(b_0,b_1) \ge R/ \sqrt{k}$, then
  \[
    \bigl| (s_0 , s_1)_{L^2} \bigr| \le  O(\ep) \, .
  \]
\end{prop}

\begin{thm} \label{in-prod}
  Let $b = b_0$ or $b_1$ and assume that
  $d( b_0, b) \le R/ \sqrt k$.
  Then 
  \begin{multline*}
    ( s_j, s_0)_{L^2} = 
    \frac 1{\sqrt{V_j V_0}}\\
    \times \int_{T^n} \exp
    \left( - \sqrt{-1} k \xt by - \frac k2 \xt(\Omega b)G(\bar{\Omega} b)
    - \frac {\sqrt{-1}}2 k \xt b (\mathrm{Re}\, \Omega)b \right) dvol\\
    + O(\ep) \, ,
  \end{multline*}
  here we take an action-angle coordinate around $\pi^{-1}(b_0)$
  as above.
  In particular, 
  \[
    \| s_i \|^2_{L^2} = 1 + O(\ep) \, .
  \]
\end{thm}

\begin{cor}
  Assume that the complex structure is invariant under the $T^n$-action
  along fibers (i.e. $G$ is constant along each fiber), then
  \[
    (s_i, s_j)_{L^2} = {\delta}_{ij} + O(\ep) \, ,
  \]
  i.e. $s_i$'s are ``asymptotically orthonormal''.
\end{cor}

\begin{rem}
  We cannot conclude that $s_i$'s are linearly independent,
  because $\dim H^0(X,L^k) = O(k^n)$ is much bigger compared to
  the error term.
\end{rem}

The next theorem is a weaker version of Theorem \ref{expansion}.
\begin{thm}\label{bs-free}
 There exists $k_0$ such that
 \[
   f_k(z) := \sum_i |s_i(s)|^2_h
   = k^n \bigl( 1+O(\ep) \bigr)
 \]
 for $k \ge k_0$.
 In particular, $s_i$'s are base point free.
\end{thm}

\subsection{Proofs}

\paragraph{Proof of Proposition \ref{decay}.}
1. 
From the choice of $r$ and $k$, there exists $y_0 \in \pi^{-1}(b_j)$
such that 
\[
    \begin{cases}
      d(z,y_0) = d \bigl( z, \pi^{-1}(b_j) \bigr) \quad \text{and}\\
      d(z,y) \ge \frac 12 d(y_0, y) 
      \quad \text{for all $y \in \pi^{-1}(b_j)$.}
    \end{cases}
\]
Hence
\begin{align*}
    d(z,y) &\ge \frac 14 \bigl( d(z,y_0) + d(y_0, y) \bigr)\\
           &= \frac 14 \bigl( d \bigl( z, \pi^{-1}(b_j) \bigr) 
              + d(y_0, y) \bigr)\, .
\end{align*}
From Theorem \ref{Th-LB},
\begin{align*}
  |s_j(z)|_h &\le \left( \frac k{2 \pi} \right)^{- \frac n4}
      \int_{\pi^{-1}(b_j)}
      |\Pi_k(z,y)| |\sigma_j(y)| dvol_y \\
    &\le \left( \frac k{2 \pi} \right)^{- \frac n4}
      \int_{\pi^{-1}(b_j)}
      Ck^n e^{-c \sqrt{k} d(z,y)} \frac 1{\sqrt{V_j}} dvol_y \\
    &\le C k^{\frac {3n}4} \frac 1{\sqrt{V_j}} 
      e^{- \frac c4 \sqrt{k} d \bigl( z, \pi^{-1}(b_j) \bigr)}
      \int_{\pi^{-1}(b_j)}
      e^{- \frac c4 \sqrt{k} d(y,y_0)} dvol_y \,.
\end{align*}
Here
\begin{align*}
  \int_{\pi^{-1}(b_j)}
  e^{- \frac c4 \sqrt{k} d(y,y_0)} dvol_y
    &\le \int_{\mathbb{R}^n} e^{-c |y|} k^{- \frac n2} dy \\
    &\le C k^{- \frac n2}.
\end{align*}
Then we obtain the desired result.

2. Since $d(z,\pi^{-1}(b_j)) \ge \frac 1{\sqrt{k}}
\left( R + \frac{3n}{4c} \log k \right)$,
\[
  e^{-c \sqrt{k} d(z,\pi^{-1}(b_j))} \le 
  e^{-c \left( R + \frac{3n}{4c} \log k \right)}
  = k^{-\frac {3n}4} e^{-cR}
  \le k^{-\frac {3n}4} \ep\, ,
\]
we have
\begin{align*}
  |s_j(z)|_h &\le \left( \frac k{2 \pi} \right)^{- \frac n4}
      \int_{\pi^{-1}(b_j)}
      Ck^n e^{-c \sqrt{k} d(z,\pi^{-1}(b_j))} 
      \frac 1{\sqrt{V_j}} dvol_y \\
    &\le C k^{\frac {3n}4} \frac 1{\sqrt{V_j}}
      \int_{\pi^{-1}(b_j)}
      k^{-\frac {3n}4} \ep dvol_y\\
    &\le C \sqrt{V_j} \ep
\end{align*}
\qed      

\paragraph{Proof of Theorem \ref{behavior}.}
To prove Theorem \ref{behavior}, we need to write down $\sigma_j$ explitely.

\begin{lem}
For $b=b_j \in \frac 1k \mathbb{Z}^n$ with $|b| \le \ep$,
  \[
    \sigma_j(y) = \frac 1{\sqrt{V_j}} \exp
    \left( - \frac{\sqrt{-1}}2 k{}^tby - 
    \frac{\sqrt{-1}}2 k{}^tb (\mathrm{Re}\, \Omega)b
    + k O(|y|^3) \right).
  \]
\end{lem}

\begin{proof}
From the choice of local trivialization of $L$,
the Hermitian metric $h(z)$ has the form
\[
  h(z) = \exp \bigl( -{}^t zG \bar{z} +O(|z|^3) \bigr)\, .
\]
So the connection on $L^k$ is
\[
  \nabla = d - k {}^t \bar{z} G dz + O(|z|^2)dz \, ,
\]
in particular, on the fiber $\pi^{-1}(b)$,
\[
  \nabla = d - k {}^t (y + \bar{\Omega} b) G dy +k O(|y|^2)dy \, .
\]
From this, the solutions of $\nabla \sigma = 0$ is of the form
\[
  \sigma = C \exp \left( \frac k2 {}^t 
  (y + \bar{\Omega} b)G(y+ \bar{\Omega}b) + k O(|y|^3) \right) e_L(y)
\]
for some constant $C$.
This constant is determined by the normalization condition 
$| \sigma_j |^2_h \equiv \frac 1{V_j}$.
\end{proof}

\paragraph{Proof of Theorem \ref{behavior}}
We put
\begin{align*}
  s_j(z) &= \left( \frac k{2 \pi} \right)^{- \frac n4}
            \int_{\pi^{-1}(b_j)} \Pi_k(z,y) \sigma_j(y) dvol(y)\\
         &= \left( \frac k{2 \pi} \right)^{- \frac n4}
            \int_{\pi^{-1}(b_j) \cap \{ |z_0 - y| \le R/2 \sqrt{k} \}}
          + \left( \frac k{2 \pi} \right)^{- \frac n4}
            \int_{\pi^{-1}(b_j) \cap \{ |z_0 - y| \ge R/2 \sqrt{k} \}}\\
         &= \left( \frac k{2 \pi} \right)^{- \frac n4}
             \bigl( (I) + (II) \bigr) \, .
\end{align*}
First we consider the second term.
\begin{align*}
  |(II)| &\le \int_{|z_0 - y| \ge R/2 \sqrt{k}}
            | \Pi_k(z,y) | \cdot |\sigma_j(y)| dvol_y\\
       &\le \int_{|z_0 - y| \ge R/2 \sqrt{k}}
            C k^n e^{-k \cdot d(z,y)^s} \frac 1{\sqrt{V_j}} dvol_y\\
       &\le C k^n \frac 1{\sqrt{V_j}} 
            \int_{\{ y \in \mathbb{R}^n \,; \,|y| \ge R/2 \} } e^{-|y|^2}
            k^{- \frac n2} dy\\
       &\le C k^{\frac n2} \frac 1{\sqrt{V_j}} e^{-R^2/4}\\
       &\le C k^{\frac n2} \frac 1{\sqrt{V_j}} \ep \, .
\end{align*}

Next we compute
\[
  (I) = \int_{d(z,y) \le R/ \sqrt k} \Pi_k (z,y) \sigma_j(y)
        dvol_y \, .
\]
Note that $dvol_y = \sqrt{\det G} dy^1 \wedge \dots \wedge dy^n$
(up to $O(1/\sqrt k)$ ).
We put $z = w/{\sqrt k}$, $b = a/\sqrt k$,
and $y = u/\sqrt k + \Omega a/\sqrt k$,
where $w \in \mathbb{C}^n$, $u \in \mathbb{R}^n$.
Then
{\allowdisplaybreaks%
\begin{align*}
  (I) &= \int_{|u| \le R}
         \Pi_k \left( \frac w{\sqrt k}, 
         \frac {u + \Omega}{\sqrt k} \right)
         \sigma_j \left( \frac {u + \Omega}{\sqrt k} \right) dvol\\
      &= \int_{|u| \le R} 
         \left( \frac k{2 \pi} \right)^n \\*
      & \phantom{ \int_{|y| \le R/} }
         \cdot \exp
         \left( - \frac 12 {}^t wG \bar{w} 
         -\frac 12 {}^t (v + \Omega a)G(v+ \bar{\Omega}a)
         + {}^t wG (v + \bar{\Omega}a) \right) \\*
      & \phantom{ \int_{|y| \le R} \exp}
         \cdot \frac 1{\sqrt{V_j}}
         \exp \left( - \frac{\sqrt{-1}}2 \xt ya - 
         \frac{\sqrt{-1}}2 \xt a (\mathrm{Re}\, \Omega)a \right) \\*
      & \phantom{ \int_{|y| \le R} \exp \exp}
         \cdot \left( 1+ O \left( \frac 1{\sqrt k} \right) \right)
         k^{-\frac n2} \sqrt{\det G} du^1 \wedge \dots \wedge du^n\\
      &= \left( \frac 1{2 \pi} \right)^n k^{\frac n2}
         \frac 1{\sqrt{V_j}} 
         \exp \left(- \frac 12 {}^t wG \bar{w} 
         -\frac {\sqrt{-1}}2 {}^t a (\mathrm{Re}\, \Omega)a \right)
         \cdot \sqrt{\det G(z)}\\*
      & \phantom{==} \cdot \int_{|u| \le R} \exp
         \left( - \frac 12 {}^t (u+ \Omega a)G(u+ \bar{\Omega} a)
         + {}^t wG(u+ \bar{\Omega} a)  
         - \frac{\sqrt{-1}}2 {}^t au \right) du\\*
      & \phantom{====}
        \cdot \left( 1+ O \left( \frac 1{\sqrt k} \right) \right).
\end{align*}
Hereafter we omit error terms.
\begin{align*}
  (I) &= \left( \frac 1{2 \pi} \right)^n 
         k^{\frac n2} \frac 1{\sqrt{V_j}} \sqrt{\det G(z)}\\*
      & \phantom{ = \frac 1{2 \pi}}
         \cdot \exp \left( - \frac 12 {}^t wG \bar{w} 
         -\frac{\sqrt{-1}}2 {}^t a (\mathrm{Re}\, \Omega)a
         - \frac 12 {}^t(\Omega a)G (\bar{\Omega} a) + 
         {}^twG(\bar{\Omega} a) \right) \\*
      & \phantom{=====} \cdot \int_{|u| \le R} \exp
         \left( - \frac 12 \xt uGu 
         - \xt vG (\mathrm{Re} \, \Omega)a + \xt uGw 
         - \frac{\sqrt{-1}}2 {}^t ua \right) du\\
      &= \left( \frac 1{2 \pi} \right)^n 
         k^{\frac n2} \frac 1{\sqrt{V_j}} \sqrt{\det G(z)}\\*
      & \phantom{\left( \frac 1{2 \pi} \right)^n}
         \cdot \exp \left( - \frac 12 \xt wG \bar{w} 
         -\frac{\sqrt{-1}}2 \xt a (\mathrm{Re}\, \Omega)a
         - \frac 12 \xt(\Omega a)G (\bar{\Omega} a) + 
         \xt wG(\bar{\Omega} a) \right) \\*
      & \phantom{\frac {k^{n/2}}{\pi^n} \frac 1{\sqrt{V_j}} \exp
         \frac {k^{n/2}}{\pi^n} \frac 1{\sqrt{V_j}}}
         \cdot \int_{|u| \le R} \exp
         \left( - \frac 12 \xt uGu + \xt uG (w-\Omega a) \right) du \, ,
\end{align*}
here we used $ \im \Omega = \frac 12 G^{-1}$.
}

\begin{lem} \label{error}
\[
  \left| \int_{\{u \in \mathbb{R}^n \, ; \,|u| \ge R\}} \exp \left(
  - \frac 12 \xt uGu +\xt vG (w-\Omega a) \right) du \right|
  \le O( \ep ) \, .
\]
\end{lem}

\begin{proof}
This Lemma follows from 
\[
  \int_{|u| \ge R} e^{-|u|^2} du \le \pi^{\frac n2} e^{- \frac n2 R^2}
\]
and the choice of $R$.
\end{proof}

\begin{lem} \label{int}
 For $\lambda \in \mathbb C$,
 \[
    \int_{\mathbb R} e^{- \frac 12 x^2 + \lambda x} dx 
    = \sqrt{2\pi} e^{\frac 12 u^2} \, .
 \]
\end{lem}
This is an easy exercise.

By combining Lemma \ref{error} and \ref{int}, we have
\begin{align*}
  &\int_{|u| \le R} \exp \left( -\frac 12 \xt uGu -\xt uG
    (\Omega a - w) \right) du \\
  & \phantom{\int_{|u| \le R}}
    = \int_{\mathbb{R}^n} \exp \left( -\frac 12 \xt uGu -\xt uG
      (\Omega a - w) \right) du + O(\ep)\\
  & \phantom{\int_{|u| \le R}}
    = \frac 1{\sqrt{det G}} (2\pi)^{\frac n2}
      \exp \left( \frac 12 \xt (\Omega a - w) G (\Omega a - w)
      \right) + O(\ep) \, .
\end{align*}

Hence, up to error terms,
{\allowdisplaybreaks
\begin{align*}
  s_j(z) &= \left( \frac 1{2 \pi} \right)^{\frac {3n}4} 
         k^{\frac n4} \frac 1{\sqrt{V_j}} \sqrt{\det G(z)}\\*
      & \phantom{\left( \frac 1{2 \pi} \right)^{\frac {3n}4}}
         \cdot \exp \left( - \frac 12 {}^t wG \bar{w} 
         -\frac{\sqrt{-1}}2 {}^t a (\mathrm{Re}\, \Omega)a
         - \frac 12 {}^t(\Omega a)G (\bar{\Omega} a) + 
         {}^twG(\bar{\Omega} a) \right) \\*
      & \phantom{\frac {k^{n/2}}{\pi^n} \frac 1{\sqrt{V_j}} \exp
         \frac {k^{n/2}}{\pi^n} \frac 1{\sqrt{V_j}}}
         \cdot \int_{\mathbb{R}^n} \exp
         \left( - \frac 12 {}^t vGv +{}^t vG (w-\Omega a) \right) dy\\
      &= \left( \frac 1{2 \pi} \right)^{\frac {3n}4} 
         k^{\frac n4} \frac 1{\sqrt{V_j}} \sqrt{\det G(z)}\\*
      & \phantom{\left( \frac 1{2 \pi} \right)^{\frac {3n}4}}
         \cdot \exp \left( - \frac 12 {}^t wG \bar{w} 
         -\frac{\sqrt{-1}}2 {}^t a (\mathrm{Re}\, \Omega)a
         - \frac 12 {}^t(\Omega a)G (\bar{\Omega} a) + 
         {}^twG(\bar{\Omega} a) \right) \\*
      & \phantom{\frac {k^{n/2}}{\pi^n} \frac 1{\sqrt{V_j}} \exp
         \frac {k^{n/2}}{\pi^n} \frac 1{\sqrt{V_j}}}
         \cdot \frac 1{\sqrt{\det G}} (\sqrt{2\pi})^n \exp
         \left( \frac 12 {}^t (w-\Omega a)G(w-\Omega a) \right)\\
      &= \left( \frac k{2\pi} \right)^{\frac n4} \frac 1{\sqrt{V_j}}
         \exp \left( - \frac 12 {}^t wG \bar{w} 
         -\frac {\sqrt{-1}}2 {}^t a (\mathrm{Re}\, \Omega)a
         - \frac 12 {}^t(\Omega a)G (\bar{\Omega} a)  \right.\\*
      & \phantom{\left( \frac {2k}{\pi} \right)^{\frac n2} 
         \frac 1{\sqrt{V_j}}  \exp \left( 
         - \frac 12 {}^t wG \bar{w} \right)}  \left.
         + {}^twG(\bar{\Omega} a) + 
         \frac 12 {}^t (w-\Omega a)G(w-\Omega a) \right)\\
      &= \left( \frac k{2\pi} \right)^{\frac n4} \frac 1{\sqrt{V_j}}
         \exp \left( \sqrt{-1} {}^t wG (\mathrm{Im}\, w) 
         - \sqrt{-1} {}^t wa - \frac 14 {}^t a G^{-1}a \right).
\end{align*}}
\qed

\paragraph{proof of Proposition \ref{inner-prod-1}.}
Let
\begin{align*}
  U_1 &= \Biggl\{ z \in X  \Biggm| 
       d(z, \pi^{-1}(b_0)) \le \frac R{2 \sqrt{k}} \Biggr\}\, ,\\
  U_2 &= \Biggl\{ z \in X  \Biggm| 
       \frac R{2 \sqrt{k}} \le d(z, \pi^{-1}(b_0))
       \le \frac 1{\sqrt{k}} \left( R + \frac{3n}{4c} \log k \right) 
       \Biggr\}\, ,\\
  U_3 &= \Biggl\{ z \in X  \Biggm| 
       d(z, \pi^{-1}(b_1)) \le \frac R{2 \sqrt{k}} \Biggr\} \, ,\\
  U_4 &= \Biggl\{ z \in X   \Biggm| 
       \frac R{2 \sqrt{k}} \le d(z, \pi^{-1}(b_1))
       \le \frac 1{\sqrt{k}} \left( R + \frac{3n}{4c} \log k \right) 
       \Biggr\}\, ,\\
  U_5 &= \Biggl\{ z \in X   \Biggm| 
       d(z, \pi^{-1}(b_0)), \, d(z, \pi^{-1}(b_1))
       \ge \frac 1{\sqrt{k}} \left( R + \frac{3n}{4c} \log k \right) 
       \Biggr\}\, .
\end{align*}
Note that $U_1 \cap U_3 = \emptyset$.
Then
\begin{align*}
  \bigl| (s_0 , s_1)_{L^2} \bigr|
  & \le \int_X \bigl| (s_0 , s_1)_h \bigr| \frac{\omega^n}{n!} \\
  & \le \int_{U_1} + \int_{U_1}
      + \int_{U_3} + \int_{U_4} + \int_{U_5} \, .
\end{align*}
From Proposition \ref{decay}, we have
\[
  \int_{U_5} \bigl| (s_0 , s_1)_h \bigr| \frac{\omega^n}{n!}
  \le \int_X O(\ep) = O(\ep)\, .
\]

From the condition $d(b_0, b_1) \ge R/ \sqrt{k}$,
if $d(z, \pi^{-1}(b_0)) \le R/ \sqrt{k}$, then
$d(z, \pi^{-1}(b_1)) \le R/ \sqrt{k}$.
Therefore, from Proposition \ref{decay},
\[
  |s_1(z)|_h \le C k^{\frac n4} 
  e^{- \frac c4 \sqrt k d(z, \pi^{-1}(b_1))}
  \le k^{\frac n4} O(\ep)\, .
\]
Hence
\begin{align*}
  \int_{U_1} \bigl| (s_0 , s_1)_h \bigr| \frac{\omega^n}{n!}
  & \le \int_{U_1} C \ep k^{\frac n2} 
     e^{- \frac c4 d(z, \pi^{-1}(b_0))} \\
  & \le C \ep k^{\frac n2} \int_{x \in \pi (U_1)}
     e^{- \frac c4 \sqrt k d(x,b_0)} \\
  & \le C \ep k^{\frac n2} 
     \int_{\bigl\{x \in \mathbb{R}^n \bigm|  |x| \le R/2 \bigr\}}
     e^{- \frac c4 |x|} k^{- \frac n2} dx \\
  & \le O(\ep)\, .
\end{align*}
Similarly we have
\[
  \int_{U_3} \bigl| (s_0 , s_1)_h \bigr| \frac{\omega^n}{n!}
  \le O(\ep)\, .
\]

Finally we consider the part of $U_2$ (and $U_4$).
If $z \in U_2 \backslash U_4$, then $|s_1(z)|_h \le O(\ep)$ by
Proposition \ref{decay}.
If $z \in U_2 \cap U_4$, then 
$|s_1(z)|_h \le k^{\frac n4} O(\ep)$
since $d(z, \pi^{-1}(b_1)) \ge R/2 \sqrt{k}$.
Consequently,
\begin{align*}
  \int_{U_2} \bigl| (s_0 , s_1)_h \bigr| \frac{\omega^n}{n!}
  & \le \int_{U_2} C \ep k^{\frac n2} 
        e^{- \frac c4 \sqrt k d(x,b_0)} \\
  & \le C \ep k^{\frac n2} \int_{\mathbb{R}^n}
        e^{- c |x|} k^{- \frac n2} dx \\
  & = O(\ep)\, .
\end{align*}
\qed

\paragraph{Proof of Theorem \ref{in-prod}.}

We put
\begin{align*}
  ( s_j, s_0 )_{L^2}
  &= \int_X \bigl( s_j(z), s_i(z) \bigr)_h \frac{\omega^n}{n!}\\
  &= \int_{d(z, \pi^{-1}(b_0)) \le \frac R{\sqrt{k}}} +
     \int_{ \frac R{\sqrt{k}} \le d(z, \pi^{-1}(b_0))
           \le \frac 1{\sqrt{k}} \left( 2R + \frac{3n}{4c} \log k\right)}\\
  & \phantom{=}
    + \int_{d(z, \pi^{-1}(b_0)) \ge 
        \frac 1{\sqrt{k}} \left( 2R + \frac{3n}{4c} \log k\right)}\\
  & = I + II + III\, .
\end{align*}

First we consider the part $III$.
Since $d(b_0,b_1) \le R/ \sqrt{k}$, if
$ d(z, \pi^{-1}(b_0)) \ge 
        \frac 1{\sqrt{k}} \left( 2R + \frac{3n}{4c} \log k\right)$,
then $ d(z, \pi^{-1}(b_1)) \ge 
        \frac 1{\sqrt{k}} \left( R + \frac{3n}{4c} \log k\right)$.
From Proposition \ref{decay},
$|s_0(z)|_h , \, |s_1(z)|_h = O(\ep)$.
Therefore,
\[
  III \le O(\ep) \, .
\]

Next we estimate $II$.
If $ d(z, \pi^{-1}(b_0)) \le 
 \frac 1{\sqrt{k}} \left( 2R + \frac{3n}{4c} \log k\right)$,
then $ d(z, \pi^{-1}(b_1)) \le 
 \frac 1{\sqrt{k}} \left( 3R + \frac{3n}{4c} \log k\right)$
In this case, we have
\[
  |z_j(z)|_h \le C k^{\frac n4} e^{-c \sqrt{k} d(z, \pi^{-1}(b_j))}.
\]
Since $d(z, \pi^{-1}(b_0)) \ge R/ \sqrt{k}$,
\[
  |s_0(z)|_h \le C k^{\frac n4} e^{-cR} \le C k^{\frac n4} \ep\, .
\]
Hence
\begin{align*}
  II &\le \int_{d(z, \pi^{-1}(b_0) \le 
      \frac 1{\sqrt{k}} \left( R + \frac{3n}{4c} \log k\right)}
      C k^{\frac n2} \ep e^{-c \sqrt{k} d(z, \pi^{-1}(b_j))}\\
     & \le C k^{\frac n2} \ep
       \int_{\mathbb{R}^n} e^{-c |u|} k^{-\frac n2} du\\
     & \le O(\ep)\, .
\end{align*}

Finally, we consider $I$.
From Theorem \ref{behavior},
\begin{align*}
  &\bigl( s_j(z), s_0(z) \bigr)_h \\
     & \phantom{s_j(z)} 
     = \left( \frac k{2\pi} \right)^{\frac n2} \frac 1{\sqrt {V_i V_j}}\\
     & \phantom{s_j(z) =} \cdot \Biggl\{
       \exp \left( \sqrt{-1} k \xt zG (\im z) -2k \sqrt{-1} \xt bz
       - \frac k4 \xt b G^{-1} b - \sqrt{-1}k \xt \bar{z}G(\im z) \right)
       \Biggr.\\
     & \phantom{ s_j(z) =
       \exp \left( \sqrt{-1} \xt zG (\im w) -2 \sqrt{-1} {}^t wa \right)
       \xt a G^{-1} a - \sqrt{-1} \xt \bar{w}G(\im w)}
       \Biggl.+ O(\ep) \Biggr\}\\     
     & \phantom{s_j(z)}
       = \left( \frac k{2\pi} \right)^{\frac n2} \frac 1{\sqrt {V_i V_j}}\\
     & \phantom{s_j(z) =} \cdot \Biggl\{
        \exp \left( -2k \xt (\im z)G(\im z) -2k \sqrt{-1} \xt bz
        -\frac k4 \xt b G^{-1} b \right)
       + O(\ep) \Biggr\}\, .
\end{align*}
We take an action-angle coordinate and write
$z = \Omega x + y$ around a point $z_0$ as above.
Then
\begin{align*}
  &\bigl( s_j(z), s_0(z) \bigr)_h \\
  & \phantom{s_j(z)} 
    = \left( \frac k{2\pi} \right)^{\frac n2} \frac 1{\sqrt {V_0 V_j}}\\
  & \phantom{s_j(z) =} \cdot \Biggl\{
    \exp \left( -k \xt x (\im \Omega) x - \sqrt{-1}k \xt b \Omega x
    - \sqrt{-1}k \xt by - \frac k2 \xt b (\im \Omega) b \right)
    \Biggr.\\
  & \phantom{ s_j(z) =
    \exp \left( \sqrt{-1} {}^t wG (\im w) -2 \sqrt{-1} {}^t wa \right)
    \xt a G^{-1} a - \sqrt{-1} \xt \bar{w}G(\im w)}
    \Biggl.+ O(\ep) \Biggr\}\, .    
\end{align*}
Hereafter we omit the error terms.
{\allowdisplaybreaks
\begin{align*}
  I &= \left( \frac k{2\pi} \right)^{\frac n2} \frac 1{\sqrt {V_0 V_j}}\\
    & \phantom{=} \cdot
       \int_{|x| \le \frac R{\sqrt{k}}, \,\, y \in T^n}
       \exp \left( -k \xt x (\im \Omega) x - \sqrt{-1}k \xt b \Omega x
       - \sqrt{-1}k \xt by - \frac k2 \xt b (\im \Omega) b \right)
       dx \, dy\\
    &= \left( \frac k{2\pi} \right)^{\frac n2} \frac 1{\sqrt {V_0 V_j}}\\
    & \phantom{=} \cdot
      \int_{y \in T^n} 
      \Biggl( \int_{|x| \le \frac R{\sqrt{k}}} 
      \exp \left(-k \xt x (\im \Omega) x 
      - \sqrt{-1}k \xt b \Omega x \right) dx \Biggr)\\
    & \phantom{=\int_{y \in T^n} }
      \exp 
      \left(- \sqrt{-1}k \xt by - \frac k2 \xt b (\im \Omega) b \right)
       dy \, .
\end{align*}
}
Here 
\begin{multline*}
  \int_{|x| \le \frac R{\sqrt{k}}} 
      \exp \left(-k \xt x (\im \Omega) x 
      - \sqrt{-1}k \xt b \Omega x \right) dx\\
  = \left( \frac k{2\pi} \right)^{-\frac n2}
    \left( \exp \left( -\frac k2 \xt b \Omega G \Omega b \right) 
    + O(\ep) \right).
\end{multline*}
Consequently,
\[
  I = \frac 1{\sqrt {V_0 V_j}} \int_{y \in T^n}
      \exp 
      \left(- \sqrt{-1}k \xt by - \frac k2 \xt b (\im \Omega) b 
      -\frac k2 \xt b \Omega G \Omega b \right) dy 
    +  O(\ep)\, . 
\]
By combining the identity
\[
  \im \Omega + \Omega G \Omega = \sqrt{-1} \re \Omega
  + \Omega G \bar{\Omega}\, ,
\]
we have
\begin{align*}
  (s_j, s_0)_{L^2}
  &= I + O(\ep)\\
  &= \frac 1{\sqrt {V_0 V_j}} \int_{T^n}
      \exp 
      \left(- \sqrt{-1}k \xt by - \frac {\sqrt{-1}}2 k \xt b (\re \Omega) b 
      -\frac k2 \xt b \Omega G \bar{\Omega} b \right) dy 
    +  O(\ep)\, . 
\end{align*}

In particular, if $s_j = s_0$ (i.e. $b = 0$),
then
\begin{align*}
  \| s_0 \|_{L^2}^2 
  &= \frac 1{V_0} \int_{\pi^{-1}(b_0)} dy +  O(\ep)\\
  &= 1 + O(\ep) \, .
\end{align*}
\qed

\paragraph{Proof of Theorem \ref{bs-free}.}
We may assume that $z \in B_{R/\sqrt{k}}(z_0)$.
Let
\begin{align*}
  U_1 &= \Biggl\{ x \in B  \Biggm| 
       d(x, \pi(z_0)) \le \frac R{ \sqrt{k}} \Biggr\}\, ,\\
  U_2 &= \Biggl\{ x \in B  \Biggm| 
       \frac R{ \sqrt{k}} \le d(x, \pi(z_0))
       \le \frac 1{\sqrt{k}} \left( R + \frac{3n}{4c} \log k \right) 
       \Biggr\}\, ,\\
  U_3 &= \Biggl\{ x \in B   \Biggm| 
       d(x, \pi(z_0))
       \ge \frac 1{\sqrt{k}} \left( R + \frac{3n}{4c} \log k \right) 
       \Biggr\}
\end{align*}
and write $f_k$ as
\begin{align*}
  f_k(z) &= \sum_{b_j \in U_1} |s_j(z)|^2_h
           + \sum_{b_j \in U_2} |s_j(z)|^2_h
           + \sum_{b_j \in U_3} |s_j(z)|^2_h\\
         &= (I) + (II) + (III)\, .
\end{align*}
From Proposition \ref{decay} we have
\[
  (III) \le \dim H^0(X,L^k) \cdot O(\ep)
        \le k^n O(\ep)
\] 
and
\begin{align*}
  (II) &\le \sum_{b_j \in U_2} 
     Ck^{\frac n2} e^{- \frac c2 \sqrt{k} d(z,\pi^{-1}(b_j))}\\
  &\le C'k^{n}
  \sum_{\genfrac{}{}{0pt}{}{b \in \frac 1k \mathbb Z^n,}
  {|b| \ge R/\sqrt{k}}}
  k^{- \frac n2} e^{-c \sqrt{k} |b|} \, .
\end{align*}
Since 
\[
  \sum_{\genfrac{}{}{0pt}{}{b \in \frac 1k \mathbb Z^n,}
  {|b| \ge R/\sqrt{k}}}
  k^{- \frac n2} e^{-c \sqrt{k} |b|}
  \longrightarrow
  \int_{\genfrac{}{}{0pt}{}{x \in \mathbb R^n,}{|x| \ge R}}
  e^{-c |x|} = C e^{-cR} = O(\ep)
\]
as $k \to \infty$, $(II)$ can be bounded from above by
$k^n O(\ep)$.

Next we compute $(I)$.
From Theorem \ref{behavior}, we have
\[
  |s_j(z)|^2_h =
  \left( \frac k{2\pi} \right)^{\frac n2} \frac 1{V_j}
  \bigl\{ \exp \left(
  -k \xt (x-b_j) (\im \Omega) (x-b_j) \right) + O(\ep) \bigr\}
\]
and
\begin{align*}
  (I) &= \left( \frac k{2\pi} \right)^{\frac n2}
  \sum_{b_j \in U_1} \frac 1{V_j}\exp \left(
  -k \xt (x-b_j) (\im \Omega) (x-b_j) \right)
  + k^n O(\ep)\\
  &= \left( \frac 1{2\pi} \right)^{\frac n2} k^n
  \sum_{\genfrac{}{}{0pt}{}
  {b \in \frac 1{\sqrt k} \mathbb Z^n + \sqrt{k} x,}
  {|b| \le R}}
  k^{-\frac n2} \frac 1{V_j}
  \exp \left(- \xt b (\im \Omega) b \right)
  + k^n O(\ep)\, .
\end{align*}
Similarly, 
\[
  \left|
  \sum_{\genfrac{}{}{0pt}{}
  {b \in \frac 1{\sqrt k} \mathbb Z^n + \sqrt{k} x,}
  {|b| \le R}}
  k^{-\frac n2} \frac 1{V_j}
  \exp \left(- \xt b (\im \Omega) b \right)
  -
  \int_{\genfrac{}{}{0pt}{}
  {u \in \mathbb R^n}{|u| \le R}}
  \exp \left(- \xt u (\im \Omega) u \right)
  \frac {du}{\sqrt{2^n \det G}}
  \right| \le \ep
\]
for large $k$ (Note that $|\sqrt{k} x| \le R$).
On the other hand, 
\[
  \left|
  \int_{|u| \ge R}
  \exp \left(- \xt u (\im \Omega) u \right)
  \frac {du}{\sqrt{2^n \det G}}
  \right| \le \ep \, .
\]
Hence
\begin{align*}
  (I) &= \left( \frac 1{2\pi} \right)^{\frac n2} k^n
         \int_{\mathbb R^n}
         \exp \left(- \xt u (\im \Omega) u \right)
         \frac {du}{\sqrt{2^n \det G}}
         + k^n O(\ep)\\
      &= \left( \frac 1{2\pi} \right)^{\frac n2} k^n
         \int_{\mathbb R^n}
         e^{-|u|^2} du
         + k^n O(\ep)\\
      &= k^n \bigl( 1 + O(\ep) \bigr) \, ,
\end{align*}
here we use $(\im \Omega)^{-1} = 2G$.
\qed



\end{document}